\documentclass[11pt,reqno]{elsarticle}

\usepackage{placeins}

\makeatletter
\def\ps@pprintTitle{%
 \let\@oddhead\@empty
 \let\@evenhead\@empty
 \def\@oddfoot{}%
 \let\@evenfoot\@oddfoot}
\makeatother

\usepackage{amsmath,amsthm,amscd,wrapfig,amsfonts,amssymb,mathtools}
\usepackage{enumitem,comment,mathabx,mathrsfs}
\usepackage{booktabs,multirow,lscape,cool,bm}
\usepackage{url,parskip,caption}
\usepackage[top=1.0in, bottom=1.0in, left=1.0in, right=1.0in]{geometry}

\newtheorem{remark}{Remark}

\makeatletter
\g@addto@macro\normalsize{%
  \setlength\abovedisplayskip{.4em}
  \setlength\belowdisplayskip{.4em}
  \setlength\abovedisplayshortskip{.4em}
  \setlength\belowdisplayshortskip{.4em}
}

\begin{document}

\begin{frontmatter}


\title
{
Phase function methods for second order inhomogeneous linear ordinary differential equations
}

\begin{abstract}

It is well known that second order homogeneous linear ordinary differential equations with
slowly varying coefficients admit slowly varying phase functions.  This observation
underlies the Liouville-Green method and many other techniques for the asymptotic
approximation of the solutions of such equations.  It is also the  basis of a recently developed
numerical algorithm that, in many cases of interest, runs in time independent
of the magnitude of the equation's coefficients and achieves accuracy on par with that
predicted by its condition number.    Here we point out that a large class of second 
order inhomogeneous linear ordinary differential equations can be efficiently and accurately solved  
by combining phase function methods for second order  homogeneous linear ordinary differential  equations
with a variant of the adaptive Levin method for evaluating oscillatory integrals.

\end{abstract}

\author[add1,add2]{Kirill Serkh}
\author[add1]{James Bremer}
\address[add2]{Department of Computer Science, University of Toronto}
\address[add1]{Department of Mathematics, University of Toronto}


\end{frontmatter}

\begin{section}{Introduction}
\label{section:introduction}

We say that $\alpha$ is a phase function for the second order homogeneous linear 
ordinary differential equation
\begin{equation}
y''(t) + q(t) y(t) = 0, \ \ \ \ \ a < t < b,
\label{introduction:hom}
\end{equation}
provided $\alpha'$ is positive on $(a,b)$ and the functions
\begin{equation}
u(t) = \frac{\cos(\alpha(t))}{\sqrt{\alpha'(t)}}
\ \ \ \mbox{and}\ \ \
v(t) =  \frac{\sin(\alpha(t))}{\sqrt{\alpha'(t)}}
\label{introduction:uv}
\end{equation}
form a basis in its space of solutions.  It is well known that when 
$q$ is a slowly varying function, (\ref{introduction:hom}) admits slowly
varying phase functions, and that this is true even when $q$
is of large magnitude and the solutions of (\ref{introduction:hom}) are highly oscillatory
or behave as combinations of rapidly increasing and decreasing exponential functions.
This observation underlies the Liouville-Green method (see, for instance, 
Chapter~6 of \cite{Olver} or Chapter~7 of \cite{Miller}) and many other techniques for the asymptotic 
approximation of the solutions of such equations
(for example, \cite{SpiglerZeros, SpiglerPhase1, SpiglerPhase2}).

It is also the basis of the algorithm of \cite{BremerKummer}, which applies when the 
coefficient $q$ is positive.  In this case,  the solutions of (\ref{introduction:hom})
are oscillatory and  the frequency of their oscillations grows with the  magnitude of $q$.
Consequently, when they are used in this regime, the running times of  standard solvers for 
ordinary differential equations grow with the magnitude of $q$.
By contrast, in many cases of interest, the numerical algorithm of \cite{BremerKummer} is able to calculate
a nonoscillatory phase function for (\ref{introduction:hom}) with near machine precision
accuracy in time independent of the magnitude of $q$.   The companion paper
\cite{BremerRokhlin} gives a bound on the complexity
of the phase function produced by the algorithm of \cite{BremerKummer} under mild
conditions on the coefficient $q$.

In this article, we describe a solver for
second order inhomogeneous liner ordinary differential equations of the form
\begin{equation}
y''(t) + q(t) y(t) = f(t),\ \ \ \ \ a < t < b,
\label{introduction:inhom}
\end{equation}
where $f$ is a slowly varying real-valued function and $q$ is slowly varying and positive.
The algorithm of \cite{BremerKummer}  is first used to construct a nonoscillatory phase function $\alpha$ such that
the functions $u$ and $v$ defined via (\ref{introduction:uv}) form a basis in the space 
of the corresponding homogeneous problem (\ref{introduction:hom}).
A particular solution of (\ref{introduction:inhom}) is then given by the formula
\begin{equation}
z(t) = u(t) \int_a^t v(s) f(s)\ ds + v(t) \int_a^t u(s) f(s)\ ds,
\label{introduction:particular}
\end{equation}
and any solution of (\ref{introduction:inhom}) is of the form
\begin{equation}
y(t) = c_1 u(t) + c_2 v(t) + z(t)
\label{introduction:general}
\end{equation}
with $c_1$ and $c_2$ constants. 
We then use an adaptive Levin method,  similar
to that  introduced in \cite{SerkhBremerLevin}, to construct a sequence
of auxillary functions $p_1,\ldots,p_m$, defined on a partition of $[a,b]$, 
which allow us to efficiently evaluate the 
integrals appearing in (\ref{introduction:particular}).  
In many cases of interest, $\alpha$ and 
the sequence of functions $p_1,\ldots,p_m$  can be constructed
in time independent of the magnitude of $q$
and, once this has been done,
constants $c_1$ and $c_2$ such that (\ref{introduction:general})
satisfies essentially any reasonable boundary conditions 
can be readily computed. 

The Levin method, which was introduced
in \cite{Levin}, is a classical technique for evaluating integrals of the form
\begin{equation}
\int_a^b \exp(i g(s)) f(s)\ ds
\label{introduction:oscint}
\end{equation}
where $f$ is slowly varying and $g$ is real-valued and slowly varying.   It was long believed that it suffers
from ``low-frequency breakdown,'' meaning that accuracy is lost when $g'$ is of insufficient magnitude.
In \cite{SerkhBremerLevin}, it is shown that, in fact, the Levin method is applicable
regardless of the magnitude of $g'$  and this observation is used to develop 
an  adaptive integration scheme for integrals of the form (\ref{introduction:oscint}).
One of the key advantages of the resulting ``adaptive Levin method'' is its ability to efficiently
evaluate oscillatory integrals in the presence of saddle points --- that is,
locations at which the first $l \geq 1$ derivatives of $g$ vanish.

Although the method of \cite{BremerKummer} can be extended to the 
case in which (\ref{introduction:hom}) has turning points (this is done in \cite{bremerphase}),
the algorithm of this paper often encounters numerical difficulties when used in this regime.
They stem from the fact that one or both of the 
functions $u$ and $v$ must be rapidly increasing in regions in which $q$ is negative and of large magnitude,
and  this makes   (\ref{introduction:general}) a numerically unstable representation
of the desired solution of (\ref{introduction:inhom}).     
Of course, in many cases in which additional information about the desired solution is known \emph{a priori}--- 
for instance, if it is a linear combination of the particular solution and a solution $u$ of 
(\ref{introduction:hom}) whose magnitude does not become large ---  a numerically stable formula
for representing the desired solution  can be devised and 
(\ref{introduction:inhom}) can be solved to high accuracy using a slightly modified version
of the algorithm presented here.

The remainder of this paper is structured as follows.  In Section~\ref{section:phase}, we  discuss
phase functions for second order homogeneous linear ordinary differential equations and describe
the algorithm of \cite{BremerKummer} for the numerical calculation of nonoscillatory phase
functions.  Section~\ref{section:levin}
discusses the adaptive Levin method for the evaluation of oscillatory integrals
and explains how a variant of it can be used to construct the sequence of functions $p_1,\ldots,p_m$
which enable the rapid evaluation of the particular solution (\ref{introduction:particular}).
In Section~\ref{section:experiments}, we  present the results of numerical experiments designed to demonstrate
the properties of our approach to solving (\ref{introduction:inhom}).     
We close with a few brief remarks in Section~\ref{section:conclusions}.

\end{section}

\begin{section}{Phase functions for second order homogeneous linear ordinary differential equations }
\label{section:phase}

\begin{subsection}{Kummer's equation}

It can be easily verified that if $\alpha$ is a phase function for the second order homogeneous linear ordinary differential equation
(\ref{introduction:hom}) --- so  that the functions $u$ and $v$ defined in (\ref{introduction:uv}) 
form a basis in its space of solutions --- then $\alpha'$ solves
\begin{equation}
q(t) - (\alpha'(t))^2  + \frac{3}{4} \left(\frac{\alpha''(t)}{\alpha'(t)}\right)^2
- \frac{1}{2} \frac{\alpha'''(t)}{\alpha'(t)} = 0,\ \ \ \ \ a <t <b.
\label{phase:kummer}
\end{equation}
Equation~(\ref{phase:kummer}) is often  referred to as Kummer's equation, after E.~E~Kummer who studied it in
\cite{Kummer}.  Kummer's equation is strongly related to the better known Riccati equation
\begin{equation}
r'(t) + (r(t))^2 + q(t) = 0,\ \ \ \ \ a < t <b,
\label{phase:riccati}
\end{equation}
satisfied by the logarithmic derivatives of the solutions of (\ref{introduction:hom}).  In particular, 
if $\alpha$ is a phase function for (\ref{introduction:hom}), then 
\begin{equation}
r(t) = i \alpha'(t) - \frac{1}{2}\frac{\alpha''(t)}{\alpha'(t)}
\end{equation}
is a solution of (\ref{phase:riccati}).

When $q$ is positive, almost all solutions of Kummer's equation are oscillatory.
However, it was recognized long ago that there exist some solutions which can 
be asymptotically approximated by nonoscillatory functions, provided $q$ is slowly varying.  This is the basis of 
many classical asymptotic techniques, including  the Liouville-Green method.   If $q$ is strictly positive
on the interval $[a,b]$, then
\begin{equation}
u_0(t) = \frac{\cos\left(\alpha_0(t)\right)}{\sqrt{\alpha_0'(t)}}
\ \ \mbox{and}\ \ \
v_0(t) = \frac{\sin\left(\alpha_0(t)\right)}{\sqrt{\alpha_0'(t)}},
\label{phase:lgbasis}
\end{equation}
where $\alpha_0$ is defined via
\begin{equation}
\alpha_0(t) = \int_a^t \sqrt{q(s)}\ ds,
\label{phase:lgphase}
\end{equation}
are a pair of Liouville-Green approximates for (\ref{introduction:hom}).
When $q$ is slowly varying, obviously so too is the function defined via (\ref{phase:lgphase}),
and this is the case regardless of the magnitude of $q$. 
For a careful discussion of the Liouville-Green method,
including rigorous error bounds for the approximates  (\ref{phase:lgbasis}),
we refer the reader to  Chapter~6 of \cite{Olver}.  A higher order asymptotic method
which operates by iteratively refining the Liouville-Green phase $\alpha_0$ is developed in
\cite{SpiglerZeros, SpiglerPhase1, SpiglerPhase2}.

A theorem showing the existence of a nonoscillatory
phase function for (\ref{introduction:hom}) in the case in which
$q(t) = \lambda^2 q_0(t)$ with $q_0$ smooth and strictly positive on $[a,b]$ appears in \cite{BremerRokhlin}.
It applies when the function $p(x) = \widetilde{p}(t(x))$, where $p(t)$ is defined via
\begin{equation}
\widetilde{p}(t) = 
\frac{1}{q_0(t)} \left(
\frac{5}{4}\left(\frac{q_0'(t)}{q_0(t)}\right)^2 - \frac{q_0''(t)}{q_0(t)}
\right)
=
4 \left(q_0(t)\right)^{\frac{1}{4}}
\frac{d}{dt} \left(
\frac{1}{\left(q_0(t)\right)^{\frac{1}{4}}}
\right)
\end{equation}
and $t(x)$ is the inverse function of 
\begin{equation}
x(t) = \int_a^t \sqrt{q_0(s)}\ ds,
\label{introduction:xt}
\end{equation}
has a rapidly decaying Fourier transform.  More explicitly, the theorem
asserts that if the Fourier transform of $p$ satisfies a bound of the form
\begin{equation}
\left|\widehat{p}(\xi)\right|\leq \Gamma \exp\left(-\mu\left|\xi\right|\right),
\end{equation}
then there exist functions $\nu$ and $\delta$ such that 
\begin{equation}
\left|\nu(t)\right| \leq  \frac{\Gamma}{2\mu} \left(1 + \frac{4\Gamma}{\lambda}\right) \exp(-\mu \lambda),
\end{equation}
\begin{equation}
\left|\widehat{\delta}(\xi)\right| \leq \frac{\Gamma}{\lambda^2}\left(1+\frac{2\Gamma}{\lambda}\right)\exp(-\mu|\xi|)
\end{equation}
and
\begin{equation}
\alpha(t) = \lambda \sqrt{q_0(t)} \int_a^t \exp\left(\frac{\delta(u)}{2}\right)\ du
\end{equation}
is a phase function for 
\begin{equation}
y''(t) + \lambda^2 \left(q_0(t) + \frac{\nu(t)}{4\lambda^2} \right) y(t) = 0.
\end{equation}
The definition of the function $p(x)$ is ostensibly quite complicated;
however, $p(x)$ is, in fact, simply equal to twice the Schwarzian derivative of the inverse function
$t(x)$ of (\ref{introduction:xt}).
This theorem ensures that, even for relatively modest values of $\lambda$,
(\ref{introduction:hom}) admits a phase function which is nonoscillatory for the purposes
of numerical computation.  Of course, when $\lambda$ is small, the frequency of oscillation of
any  phase function for (\ref{introduction:hom}) is modest.

\vskip 1em
\begin{remark}
When discussing phase function for second order linear ordinary differential
equations it is standard practice to restrict attention to the special form 
(\ref{introduction:hom}).  This is because Kummer's equation and the phase function
$\alpha$ are invariant under the standard transformation which takes 
the more general ODE
\begin{equation}
y''(t) + p(t) y'(t) + q(t) y(t) = 0,\ \ \ \ \ a < t < b,
\label{phase:generalode}
\end{equation}
to the form  (\ref{introduction:hom}).  To see this, we suppose that
$\alpha$ is a phase function for  (\ref{phase:generalode}) so that
$\alpha'$ is positive on $(a,b)$ and 
\begin{equation}
\sqrt{\frac{\omega(t)}{\alpha'(t)}} \cos\left(\alpha(t)\right)
\ \ \ \mbox{and}\ \ \ 
\sqrt{\frac{\omega(t)}{\alpha'(t)}} \sin\left(\alpha(t)\right),
\label{phase:generaluv}
\end{equation}
where
\begin{equation}
\omega(t) = \exp\left(-\int p(t)\ dt \right),
\end{equation}
form a basis in its space of solutions.  We note that by Abel's theorem, the Wronskian
of any pair of independent solutions of (\ref{phase:generalode})
is a constant multiple of $\omega$, and so the factors of $\sqrt{\omega(t)}$ appearing
in (\ref{phase:generaluv}) are necessary.    It can be easily shown that $\alpha'$ satisfies
\begin{equation}
q(t)- \frac{(p(t))^2}{4} - \frac{p'(t)}{2} - (\alpha'(t))^2  + \frac{3}{4} \left(\frac{\alpha''(t)}{\alpha'(t)}\right)^2
- \frac{1}{2} \frac{\alpha'''(t)}{\alpha'(t)} = 0,\ \ \ \ \ a <t <b.
\label{phase:generalkummer}
\end{equation}

If we let 
\begin{equation}
z(t) = \exp\left(-\frac{1}{2}\int p(t)\ dt\right) y(t),
\end{equation}
where the choice of antiderivative is immaterial,
then $z$ solves the second order linear ordinary differential equation
\begin{equation}
z''(t) + \left(q(t) - \frac{(p(t))^2}{4} - \frac{p'(t)}{2} \right) z(t) = 0
\label{phase:ode2}
\end{equation}
and substituting the coefficient in (\ref{phase:ode2}) into Kummer's equation
(\ref{phase:kummer}) simply yields (\ref{phase:generalkummer}).
\end{remark}

\end{subsection}


\begin{subsection}{A numerical algorithm for constructing nonoscillatory phase functions}
\label{section:phase:algorithm}

In this subsection, we describe the numerical method of \cite{BremerKummer}
for constructing a nonoscillatory phase
function for equations of the form (\ref{introduction:hom}) in the event that the coefficient
$q$ is slowly varying and positive on the interval $[a,b]$.  The principal difficulty in designing such
an algorithm is that while there exists a nonoscillatory solution of (\ref{phase:kummer}),
almost all solutions of Kummer's equation are oscillatory and some mechanism must be used
to identify a nonoscillatory solution.    
The algorithm of  \cite{BremerKummer} addresses this problem  by introducing 
a ``windowed'' version of the coefficient  $\widetilde{q}$ which is
equal to a  constant $\nu^2$ on the right quarter of the interval $[a,b]$ and nearly equal to 
the original coefficient $q$ on the left quarter of $[a,b]$.   The equation
\begin{equation}
y''(t) + \widetilde{q}(t) y(t) = 0
\end{equation}
admits a nonoscillatory phase function $\widetilde{\alpha}$ which is equal to $i \nu t $ on 
the left quarter of $[a,b]$ and closely approximates a nonoscillatory phase function
for $q$ on the right quarter of $[a,b]$.  Accordingly, the algorithm solves an initial value
problem for Kummer's equation 
 to calculate the values of a nonoscillatory phase function for (\ref{introduction:hom})
and its first few derivatives at the point $b$.    A terminal value problem for Kummer's equation
is then solved to construct the desired phase function over the entire interval $[a,b]$.

We now give a more detailed description of the algorithm of \cite{BremerKummer}.
It takes as inputs the endpoints $a$ and $b$,  an external subroutine for evaluating the coefficient $q$,
an integer parameter $k$ which determines the order of the piecewise Chebyshev expansions 
used to represent the outputs and a real-valued parameter $\epsilon$ which specifies the desired precision
for the calculation.  The outputs of the algorithm comprise three $k^{th}$ order piecewise Chebyshev expansions
which represent the phase function $\alpha$ and its first two derivatives $\alpha'$ and $\alpha''$.
To be entirely clear, by an  $k^{th}$ order piecewise Chebyshev 
expansions  on the interval $[a,b]$, we mean  a sum of the form
\begin{equation}
\begin{aligned}
&\sum_{i=1}^{m-1} \chi_{\left[x_{i-1},x_{i}\right)} (t) 
\sum_{j=0}^{k} c_{ij}\ T_j\left(\frac{2}{x_{i}-x_{i-1}} t + \frac{x_{i}+x_{i-1}}{x_{i}-x_{i-1}}\right)\\
+
&\chi_{\left[x_{m-1},x_{m}\right]} (t) 
\sum_{j=0}^{k} c_{mj}\ T_j\left(\frac{2}{x_{m}-x_{m-1}} t + \frac{x_{m}+x_{m-1}}{x_{m}-x_{m-1}}\right),
\end{aligned}
\end{equation}
where $a = x_0 < x_1 < \cdots < x_m = b$ is a partition of $[a,b]$,
$\chi_I$ is the characteristic function on the interval $I$ and 
$T_j$ denotes the Chebyshev polynomial of degree $j$.

We let
\begin{equation}
\nu = \sqrt{ q\left(\frac{a+b}{2}\right)}
\end{equation}
and define  $\widetilde{q}$ via the formula
\begin{equation}
\widetilde{q}(t) = \phi(t) \nu^2 + (1-\phi(t)) q(t),
\end{equation}
where $\phi$ is given by 
\begin{equation}
\phi(t) = \frac{
1+\mbox{erf}\left(\frac{12}{b-a} \left(t-\frac{a+b}{2}\right)\right)
}{2}.
\label{algorithm:errfun}
\end{equation}
The constant in (\ref{algorithm:errfun}) is chosen so that 
\begin{equation}
\left|\phi(a) \right|,\ \left|\phi(b) - 1\right| < \epsilon_0,
\end{equation}
where $\epsilon_0$ denotes IEEE double precision machine zero.  
We next solve the terminal value problem
\begin{equation}
\left\{
\begin{aligned}
&\widetilde{q}(t) - (\widetilde{\alpha}'(t))^2 + \frac{3}{4} \left(\frac{\widetilde{\alpha}''(t)}
{\widetilde{\alpha}'(t)}\right)^2
- \frac{1}{2} \frac{\widetilde{\alpha}'''(t)}{\widetilde{\alpha}'(t)} = 0, \ \ \ \ \ a <t <b,\\
&\widetilde{\alpha}'(b) = \nu\\
&\widetilde{\alpha}''(b) = 0
\end{aligned}
\right.
\label{algorithm:windowkummer}
\end{equation}
using the algorithm described in \ref{section:algorithm:odesolver}.
We take the precision parameter for that algorithm to be the input parameter $\epsilon$ and the integer parameter 
specifying the order of the piecewise Chebyshev expansions it outputs to be the input parameter $k$.
 Although the algorithm produces $k^{th}$ order piecewise Chebyshev expansions representing
the functions  $\widetilde{\alpha}'$ and  $\widetilde{\alpha}''$, it is only the 
values of these functions at the point $a$ which concern us. 

Next, we use the method of  \ref{section:algorithm:odesolver} to solve the initial value problem
\begin{equation}
\left\{
\begin{aligned}
&q(t) - (\alpha'(t))^2 + \frac{3}{4} \left(\frac{\alpha''(t)}
{\alpha'(t)}\right)^2
- \frac{1}{2} \frac{\alpha'''(t)}{\alpha'(t)} = 0,\ \ \ \ \ a < t <b, \\
&\alpha'(a) = \widetilde{\alpha}'(a) \\
&\alpha''(a) = \widetilde{\alpha}''(a).
\end{aligned}
\right.
\label{algorithm:windowkummer2}
\end{equation}
Once again, the precision parameter is taken to be $\epsilon$ and the integer parameter 
is taken to be $k$.   The outputs of the procedure consist of $k^{th}$ order piecewise
Chebyshev expansions representing $\alpha'$ and $\alpha''$.
A $k^{th}$ order piecewise Chebyshev expansion
representing the phase function $\alpha$ itself is constructed via spectral integration;
the particular choice of antiderivative is largely irrelevant, and we 
determine it through the requirement that  $\alpha(a) = 0$.

\end{subsection}

\end{section}

\begin{section}{Levin methods}
\label{section:levin}

\begin{subsection}{The classical Levin method}

The classical Levin method \cite{Levin} is a technique for rapidly evaluating integrals
of the form
\begin{equation}
\int_a^b \exp(i g(t)) f(t)\ dt
\label{levin:oscint}
\end{equation}
in the event that $f$ is a slowly varying function, g is a slowly varying real-valued function and 
$g'$ is of  large magnitude.  It proceeds by constructing a solution to the ordinary differential 
equation
\begin{equation}
p'(t) + i g'(t) p(t) = f(t),\ \ \ \ \  a < t < b.
\label{levin:ode}
\end{equation}
Then
\begin{equation}
\frac{d}{dt} \left( \exp(ig(t)) p(t) \right) = \exp(i g(t)) f(t)
\end{equation}
and the value of (\ref{levin:oscint}) is given by
\begin{equation}
p(b)\exp(i g(b)) - p(a)\exp(i g(a)).
\label{levin:value}
\end{equation}
This procedure is more efficient than standard methods under the above assumptions because,
as shown in \cite{Levin}, (\ref{levin:ode}) admits a nonoscillatory solution $p_0$ in this case.

The differential equation (\ref{levin:ode}) is typically solved via a spectral collocation method.
Although the operator
\begin{equation}
D\left[p\right](t) = p'(t) + i g'(t) p(t)
\label{levin:op}
\end{equation}
appearing on the left-hand side of (\ref{levin:ode}) has a nontrivial nullspace, the resulting discretization
of it  will be well-conditioned as long an appropriate discretization scheme is used.   In particular, 
it is necessary  to chose a collocation grid sufficient to resolve $f$ and $g$
but not the nullspace of the operator $D$, which consists of all multiples of the function 
\begin{equation}
\eta(t) = \exp(-i g(t)).
\label{levin:eta}
\end{equation}
When $g'$ is not of sufficiently large magnitude, it is not possible to 
choose such a grid of collocation points and, in this event, the matrix discretizing (\ref{levin:ode}) will have a small singular value
and be ill-conditioned.  This phenomenon is known as ``low-frequency breakdown'' and it 
has long been viewed as a barrier to applying the Levin method in such cases.

\end{subsection}

\begin{subsection}{The adaptive Levin method}

The articles \cite{LevinLi,LiImproved} present experimental evidence indicating that when a Chebyshev spectral
method is used to discretize (\ref{levin:ode}) and  the resulting linear system is 
solved via a truncated singular value decomposition, no low-frequency breakdown seems to occur
in practice.  In \cite{SerkhBremerLevin}, a proof that this is, in fact, the case is presented
and it is observed that the lack of low-frequency breakdown makes it possible to use
the Levin method as the basis of an adaptive integration scheme.  More explicitly, because
subdividing the interval $[a,b]$ over which (\ref{levin:oscint}) is given has the 
practical effect of reducing the  magnitude of $g'$,  an adaptive algorithm is only
viable if high accuracy can be achieved even when the magnitude of $g'$ is small.  One of the remarkable features of the 
adaptive Levin method introduced in  \cite{SerkhBremerLevin} is that it is still quite efficient
when $g$ has saddle points; i.e., locations where one or more of the derivatives of $g$ vanishes.  

The procedure of \cite{SerkhBremerLevin} is quite simple --- it is completely analogous to an
adaptive Gaussian quadrature scheme, but on each subinterval $[a_0,b_0]$ considered,
it uses a Chebyshev spectral method to solve 
\begin{equation}
p'(t) + i g'(t) p(t) = f(t), \ \ \ \ \ a_0 < t < b_0,
\label{levin2:ode}
\end{equation}
and estimates the value of 
\begin{equation}
\int_{a_0}^{b_0} \exp(i g(x)) f(x)\ dx
\label{levin2:oscint}
\end{equation}
via the difference
\begin{equation}
p(b) \exp(i g(b)) -  p(a) \exp(i g(a)) 
\end{equation}
rather than using a Gaussian quadrature rule to estimate (\ref{levin2:oscint}) more directly.

\end{subsection}

\begin{subsection}{Numerical algorithm}
\label{section:levin:algorithm}

In this subsection, we describe a variant of the adaptive Levin algorithm for 
efficiently evaluating the particular solution (\ref{introduction:particular}).
It operates by  constructing a collection of  functions $p_1,\ldots,p_m$ that
allow for the efficient calculation of 
\begin{equation}
\int_a^t \exp(i g(s)) \widetilde{f}(s)\ ds,
\label{levin:oscint2}
\end{equation}
where
\begin{equation}
g(s) = \alpha(s)
\label{levin:g}
\end{equation}
and
\begin{equation}
\widetilde{f}(s) = \frac{f(s)}{\sqrt{\alpha'(s)}}.
\label{levin:ftilde}
\end{equation}
Since
\begin{equation}
\int_a^t \exp(i g(s)) \widetilde{f}(s)\ ds
=
\int_a^t \frac{\exp(i \alpha(s))}{\sqrt{\alpha'(s)}} f(s)\ ds
= \int_a^t u(s) f(s)\ ds  + i\int_a^t v(s) f(s)\ ds,
\end{equation}
the particular solution defined via (\ref{introduction:particular})
can be evaluated at a point $t$ if the values of 
$u(t)$, $v(t)$ and the integral (\ref{levin:oscint2}) are known.

The algorithm takes as input a subroutine
for evaluating the phase function $\alpha$ and its first derivative
as well as the input function $f$,
a real-valued parameter $\epsilon$ specifying the desired accuracy for the computations
and a positive integer $k$.  The output of the algorithm comprises a partition
\begin{equation}
a=a_0 < a_1 < \ldots < a_m = b
\label{levin:partition}
\end{equation}
of the interval $[a,b]$ and a collection of $m$ Chebyshev expansions.
Each Chebyshev expansion is of order $(k-1)$ and the $j^{th}$ expansion 
represents the $j^{th}$ output function $p_j$, 
which is defined over the interval $[a_{j-1},a_j]$ and closely approximates
a solution of the differential equation
\begin{equation}
p'(t) + i g'(t) p(t) = \widetilde{f}(t), \ \ \ \ \ a_{j-1} < t < a_j. 
\end{equation}
It follows that (\ref{levin:oscint2}) can be calculated
by finding the least positive integer $j$  such that $t \leq a_j$ and then computing the sum
\begin{equation}
\biggl[ p_{j}(t) \exp(i g(t)) -p_j(a_{j-1}) \exp(ig(a_{j-1}))\biggr]
 + \sum_{i=1}^{j-1} \biggl[ p_i(a_{i}) \exp(ig(a_i)) -p_i(a_{i-1}) \exp(ig(a_{i-1})) \biggr].
\label{levin:sum}
\end{equation}

The algorithm maintains a list of subintervals of $[a,b]$ which have been processed and
a list of ``accepted'' subintervals of $[a,b]$.   Upon completion of the procedure,
the list of  ``accepted subintervals'' determines the partition (\ref{levin:partition}).
Initially, the list of subintervals to process contains $[a,b]$ and the list of accepted subintervals
is empty.  The following 
sequence of operations is performed repeatedly until the 
list of subintervals to process is empty:
\begin{enumerate}

\item
Remove a subinterval $[a_0,b_0]$ from the list of subintervals to process.

\item Form the $k$-point extremal Chebyshev grid $x_1,\ldots,x_k$ on the interval $[a_0,b_0]$;
that is, let 
\begin{equation}
x_j = \frac{b_0-a_0}{2} \cos\left(\pi \frac{k-j}{k-1}\right)   + \frac{b_0+a_0}{2},\  \ \ \ \ j=1,\ldots,k.
\label{algorithm:nodes}
\end{equation}

\item
Form the $k \times k$ Chebyshev spectral differentiation matrix
$B$ which takes the vector
\begin{equation}
\left(
\begin{array}{c}
h(x_1) \\
h(x_2)\\
\vdots\\
h(x_k)
\end{array}
\right)
\end{equation}
of values of any expansion of the form
\begin{equation}
h(t) = \sum_{j=0}^{k-1} 
c_j T_j\left(
\frac{2}{b_0-a_0} t + \frac{b_0+a_0}{b_0-a_0}
\right)
\end{equation}
at the extremal Chebyshev nodes  to the vector 
\begin{equation}
\left(
\begin{array}{c}
h'(x_1) \\
h'(x_2)\\
\vdots\\
h'(x_k)
\end{array}
\right)
\end{equation}
of the values of its derivative at the extremal Chebyshev nodes.

\item
Call the external subroutine provided as input to evaluate $g$, $g'$ and $f$ at the points $x_1,\ldots,x_k$
and 
form the $k \times k$ matrix
\begin{equation}
A = B + 
i \left(\begin{array}{ccccccc}
g'(x_1)&       &        &\\
        & g'(x_2) &        &\\
        &         & \ddots \\&
        &         &          & g'(x_k)
\end{array}
\right)
\label{levin:specA}
\end{equation}
which discretizes the restriction of the differential operator
$D$ defined in (\ref{levin:op}) to the interval $[a_0,b_0]$.

\item Solve the linear system 
\begin{equation}
A \left(\begin{array}{c}y_1\\y_2\\ \vdots \\ y_k\end{array}\right)
=\left(
\begin{array}{c}
\widetilde{f}(x_1)\\
\widetilde{f}(x_2)\\
\vdots\\
\widetilde{f}(x_k)
\end{array}
\right),
\label{levin:speceq}
\end{equation}
where  $\widetilde{f}$ is defined in (\ref{levin:ftilde}),
using a truncated singular value decomposition.
To do so, first  construct the singular value decomposition
\begin{equation}
A = 
\left(
\begin{array}{cccc}
U_1 & U_2 & \cdots & U_k
\end{array}  
\right)
\left(
\begin{array}{ccccc}
\sigma_1 &           &         & \\
         & \sigma_2  &         & \\
         &           & \ddots  & \\
         &           &         &\sigma_k \\
\end{array}
\right) 
\left(
\begin{array}{cccc}
V_1 & V_2 & \cdots & V_k
\end{array}  
\right)^*
\end{equation}
of the matrix $A$.    Then, find the least integer $1\leq l \leq k$ such that $\sigma_l > 10 \epsilon_0 \|A\|_F$, where $\epsilon_0$
is machine zero and $\|A\|_F$ is the Frobenius norm of $A$.    Finally,  let
\begin{equation}
\left(
\begin{array}{c}
y_1\\
y_2\\
\vdots\\
y_k
\end{array}
\right)
  = 
\left(
\begin{array}{ccccc}
V_1 & \cdots  & V_l
\end{array}
\right)
\left(
\begin{array}{ccccc}
\sigma_1 &           &         & \\
         & \sigma_2  &         & \\
         &           & \ddots  & \\
         &           &         &\sigma_l \\
\end{array}
\right) ^{-1}
\left(\begin{array}{cccccc}
U_1 & \cdots & U_l 
\end{array}\right)^*
\left(
\begin{array}{c}
\widetilde{f}(x_1) \\
\widetilde{f}(x_2) \\
\vdots\\
\widetilde{f}(x_k)
\end{array}
\right).
\label{levin:solve}
\end{equation}

\item  
Calculate the coefficients $c_0,\ldots,c_{k-1}$
such that the Chebyshev expansion defined via
\begin{equation}
h(t) = \sum_{j=0}^{k-1} 
c_j T_j\left(
\frac{2}{b_0-a_0} t + \frac{b_0+a_0}{b_0-a_0}
\right)
\label{levin:chebexp}
\end{equation}
satisfies 
\begin{equation}
h(x_1) = y_1, \ h(x_2) = y_2,\ \ldots,\ h(x_k)=y_k.
\end{equation}

\item
If 
\begin{equation}
\frac{\sum_{j=\lfloor k/2 \rfloor }^{k-1} c_j^2 }{\sum_{j=0}^{k-1} c_j^2} < \epsilon^2,
\end{equation}
then add the interval $[a_0,b_0]$ to the list of accepted intervals.
In this case, the interval $[a_0,b_0]$ becomes one of the subintervals
$[a_j,a_{j+1}]$  in the partition (\ref{levin:partition}) and 
(\ref{levin:chebexp}) becomes the output function $p_j$.
Otherwise, add the intervals $[a_0,c_0]$ and $[c_0,b_0]$, where
$c_0 = (a_0+b_0)/2$, to the list of subintervals to process.
\end{enumerate}
\vskip 1em

When the spectral discretization on an interval $[a_0,b_0]$
suffices to resolve $f$, $g$ and the nonoscillatory solution $p_0$ whose existence is 
established in \cite{Levin}, the $(k-1)^{st}$ order Chebyshev expansion
which agrees with $p_0$ at the $k$-point Chebyshev extremal grid on $[a_0,b_0]$ 
approximates a solution of (\ref{levin:speceq}) with high accuracy.  
If, in addition, (\ref{levin:specA}) is invertible --- which is the case
when $g'$ is of sufficiently large magnitude --- it follows that 
the obtained solution (\ref{levin:chebexp})  will closely approximate
$p_0$.  So when $g'$ is of sufficient magnitude, the functions $p_1,\ldots,p_m$
produced by the procedure just described comprise a piecewise Chebyshev
discretization of the Levin solution $p_0$ of the equation
(\ref{levin:ode}).

However, when $g'$ is not of sufficient magnitude, the matrix
(\ref{levin:specA}) is noninvertible and the obtained solution
(\ref{levin:chebexp}) is a linear combination of $p_0$ and some other function
(when the nullspace is fully resolved, that other function is a multiple 
of (\ref{levin:eta}), but this need
not be the case when the discretization fails to fully resolve $\eta$).
So when $g'$ is  of insufficient magnitude, there is no guarantee that the functions $p_j$
will approximate $p_0$ or even that $p_j(a_j)$ will equal $p_{j+1}(a_j)$.
Nonetheless, Formula~(\ref{levin:sum}) remains an accurate and efficient
mechanism for evaluating (\ref{levin:oscint2}).
We note, however, that 
because $p_{j-1}(a_j)$ need not be equal to $p_j(a_j)$, it is not possible
to apply the obvious simplification to (\ref{levin:oscint2}).

\vskip 1em
\begin{remark}
The truncated singular value decomposition is quite expensive.  In our implementation
of the method described in this subsection, we used a rank-revealing QR decomposition in lieu
of the truncated singular value decomposition to 
solve the linear system which results from discretizing the Levin equation.
This was found to be about 5 times faster and lead to no apparent loss in accuracy.
\end{remark}

\end{subsection}

\end{section}

\begin{figure}[h!!!!!!!!!!!!!!!!!!]
\centering
\hfil
\includegraphics[width=.65\textwidth]{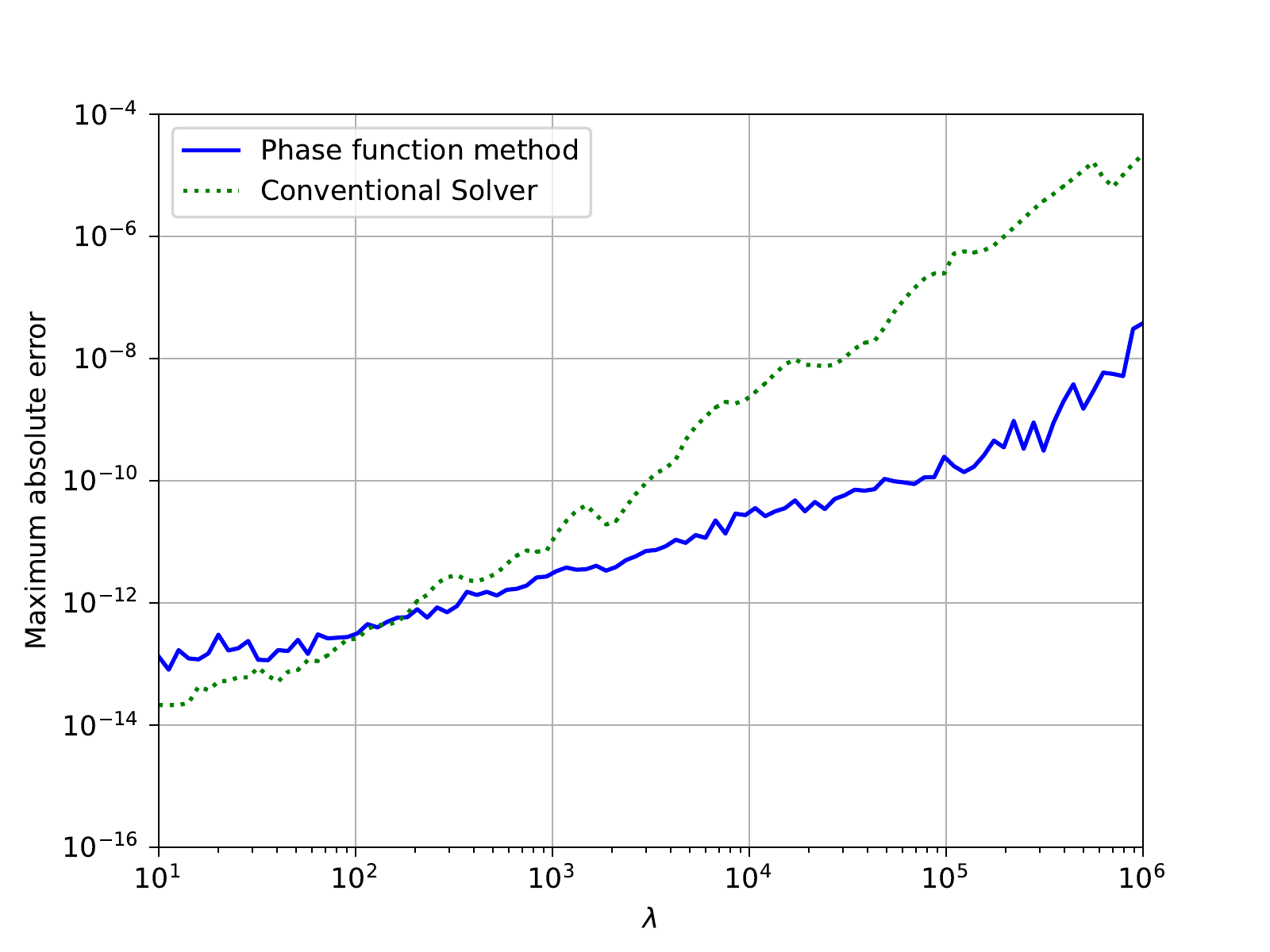}
\hskip 3em
\hfil

\hfil
\includegraphics[width=.40\textwidth]{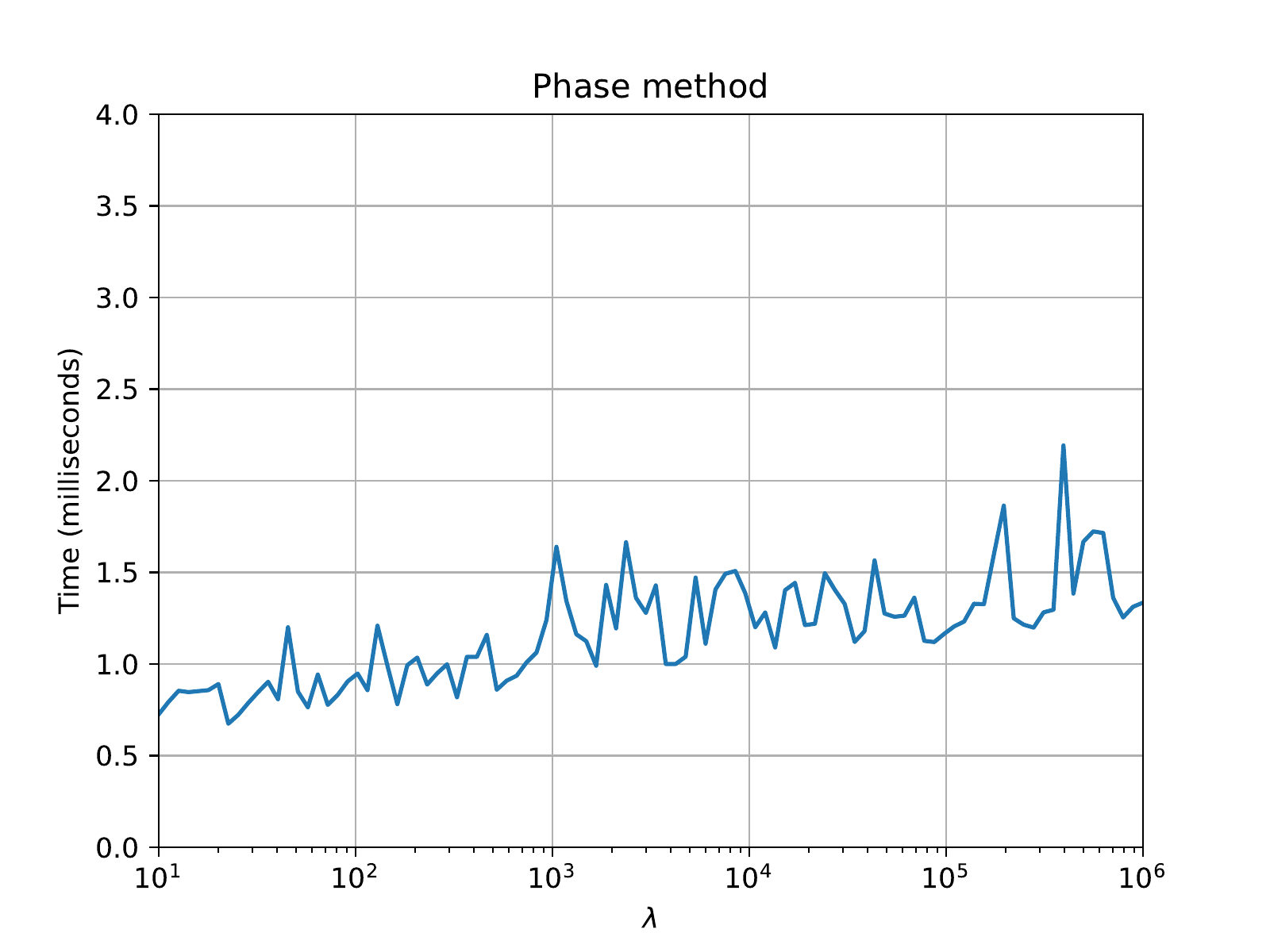}
\hfil 
\includegraphics[width=.40\textwidth]{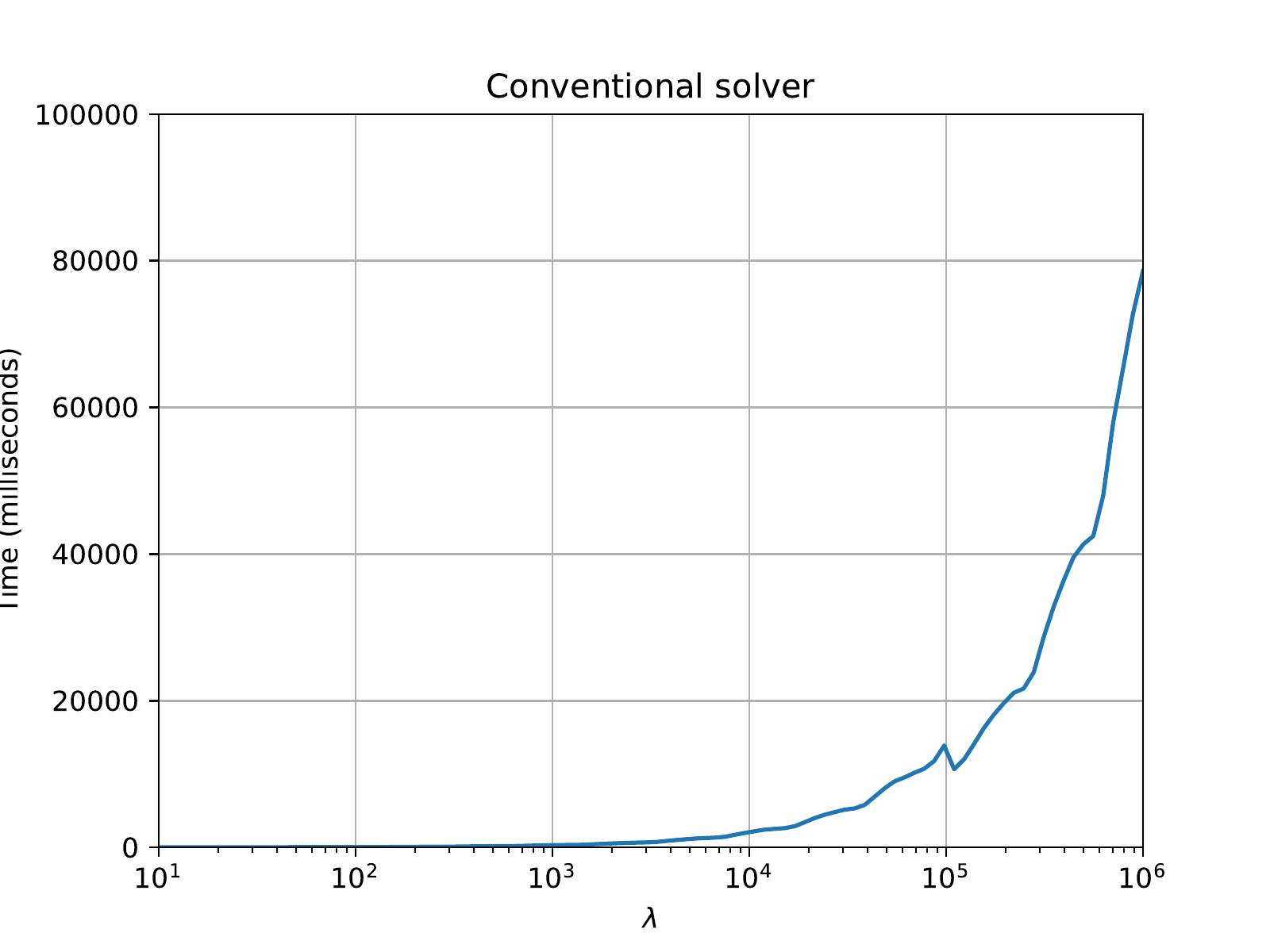}
\hfil

\hfil
\includegraphics[width=.40\textwidth]{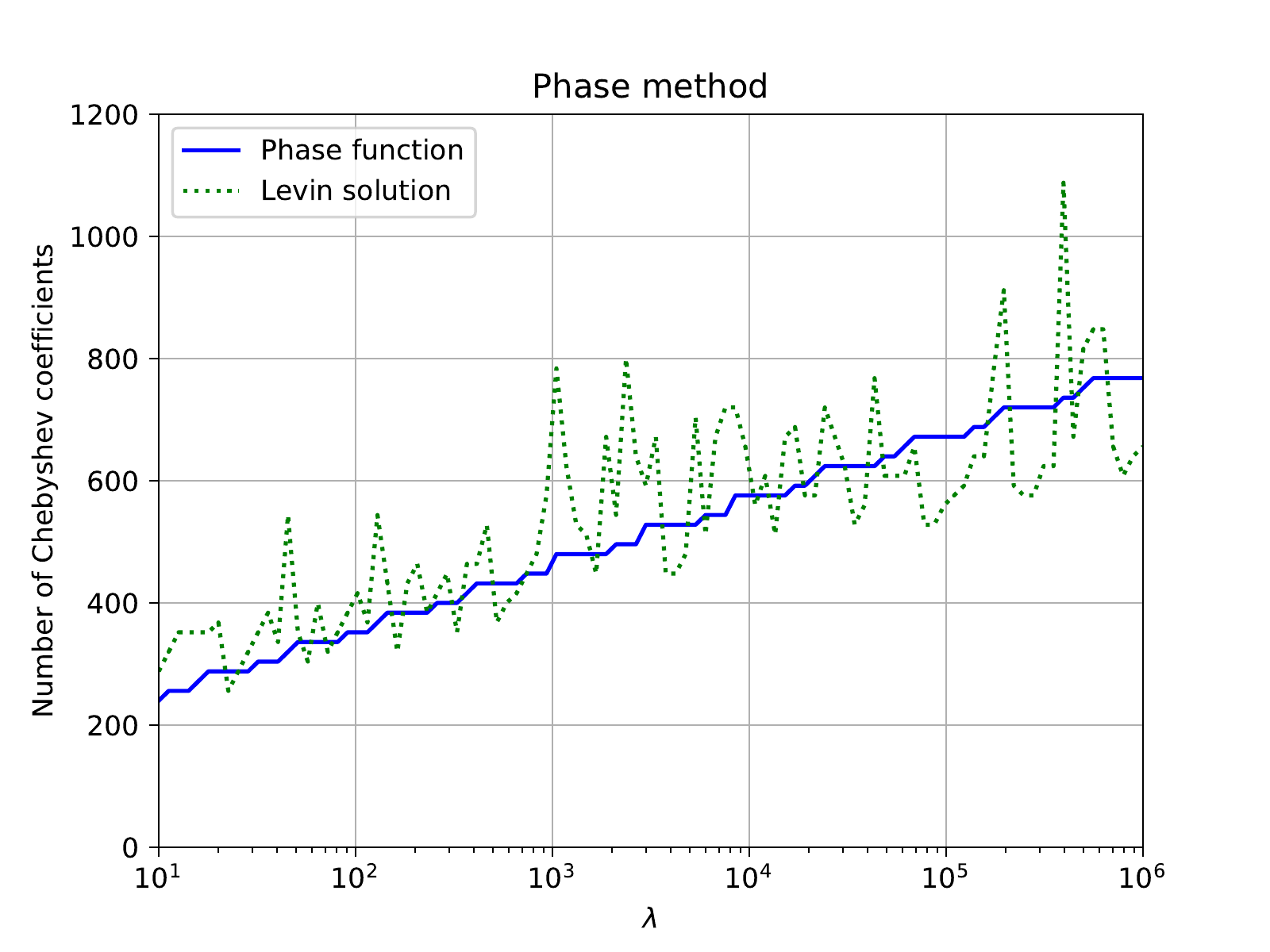}
\hfil
\includegraphics[width=.40\textwidth]{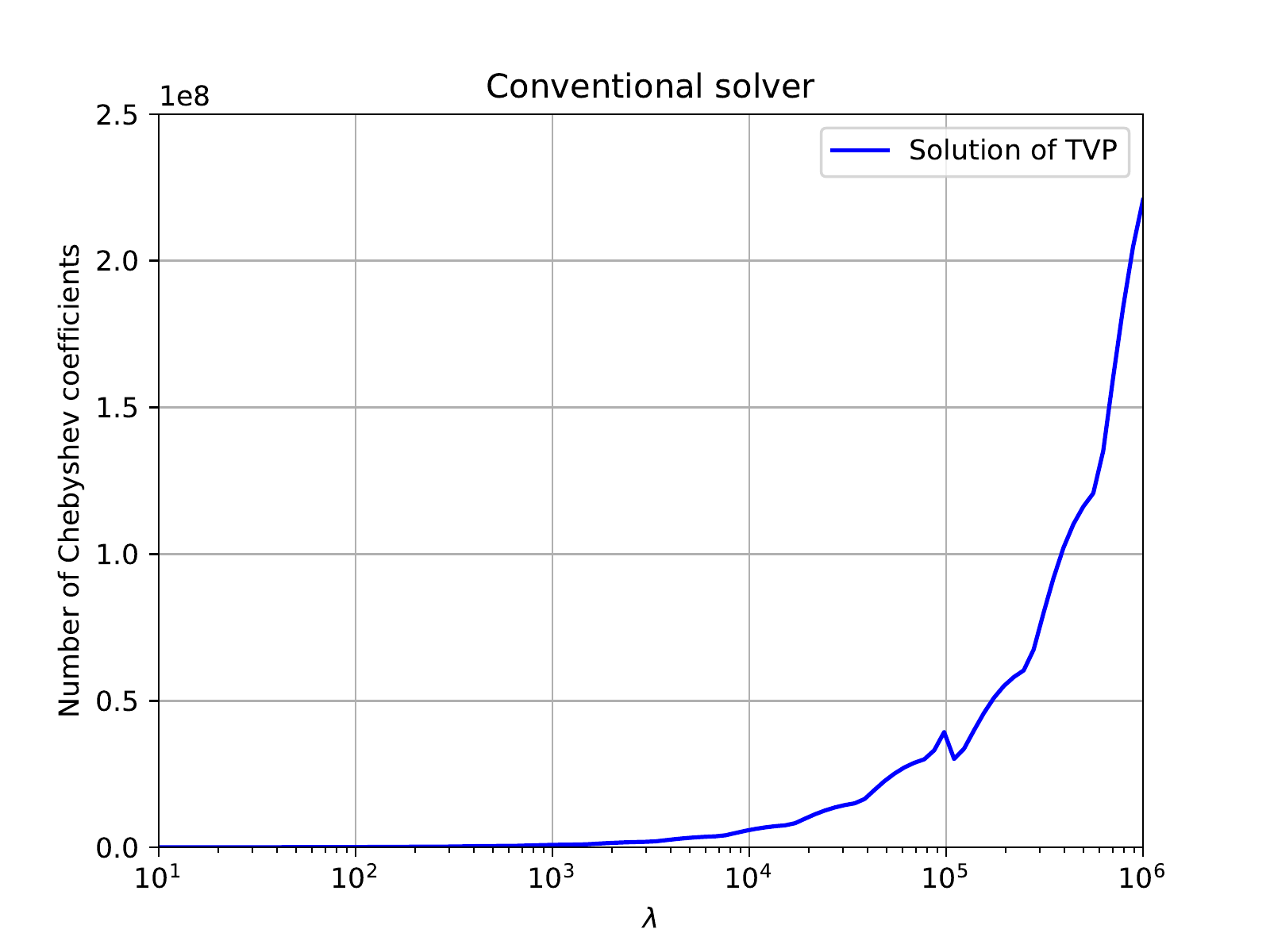}
\hfil

\caption{
The results of the first set of experiments of Subsection~\ref{section:experiments:airy}.
The plot in the first row gives the maximum observed absolute errors in the solutions
obtained by the phase function method and the conventional solver as functions of $\lambda$.
The plot on the left-hand side of the second row gives the time required by the phase method as
a function of $\lambda$, while the plot on the right-hand side of the second row gives the time
required by the conventional solver as a function of $\lambda$.
The plot on the left-hand side of the third row gives the number of coefficients
in the piecewise Chebyshev expansions of the phase function and the Levin solution
as functions of $\lambda$, whereas the plot on the right-hand side of the third row
gives the number of coefficients in the piecewise Chebyshev expansion
of the solution of (\ref{experiments:airy:ivp}) produced by the conventional solver.
}
\label{figure:airyplots1}
\end{figure}

\begin{section}{Numerical experiments}
\label{section:experiments}

In this section, we present the results of numerical experiments
which were conducted to illustrate the properties of the
algorithm of this paper.  The code for these experiments
was written in Fortran and compiled with version 12.10
of the GNU Fortran compiler.     They were performed on
a workstation computer equipped  with an AMD 3995WX  processor and 512GB
of memory.  No attempt was made to parallelize our code (i.e., only a single processor
core was used).

In the course of conducting the experiments for this article we found
that, in most cases, our  method obtains higher accuracy than the 
conventional solver described in \ref{section:algorithm:odesolver}.
Accordingly, in almost all of the experiments described here, we measured the accuracy
of the solutions obtained via our algorithm by comparison with 
reference solutions calculated by running the conventional solver of
\ref{section:algorithm:odesolver} using quadruple precision arithmetic
(i.e., using Fortran REAL*16 numbers).
These extended precision calculations required a great deal of time and memory, which
is the reason that we conducted the experiments on a workstation computer
equipped with a large amount of memory.
The experiments of Subsection~\ref{section:experiments:airy},
which concern a terminal value  problem whose solution is explicitly known,
are the exceptions.

In all of the experiments discussed here, 
$15^{th}$ order piecewise Chebyshev expansions
were used to represent phase functions as well as the functions $p_1,\ldots,p_m$, and 
the tolerance  parameter $\epsilon$ passed to the algorithms of Subsections~\ref{section:phase:algorithm}
and \ref{section:levin:algorithm} was taken to be to $10^{-13}$.
Whenever the conventional solver was deployed, $15^{th}$ order 
piecewise Chebyshev expansions were used to represent the obtained solutions.

\begin{subsection}{A terminal value problem with a known solution}
\label{section:experiments:airy}

\begin{figure}[t!]

\hfil
\includegraphics[width=.40\textwidth]{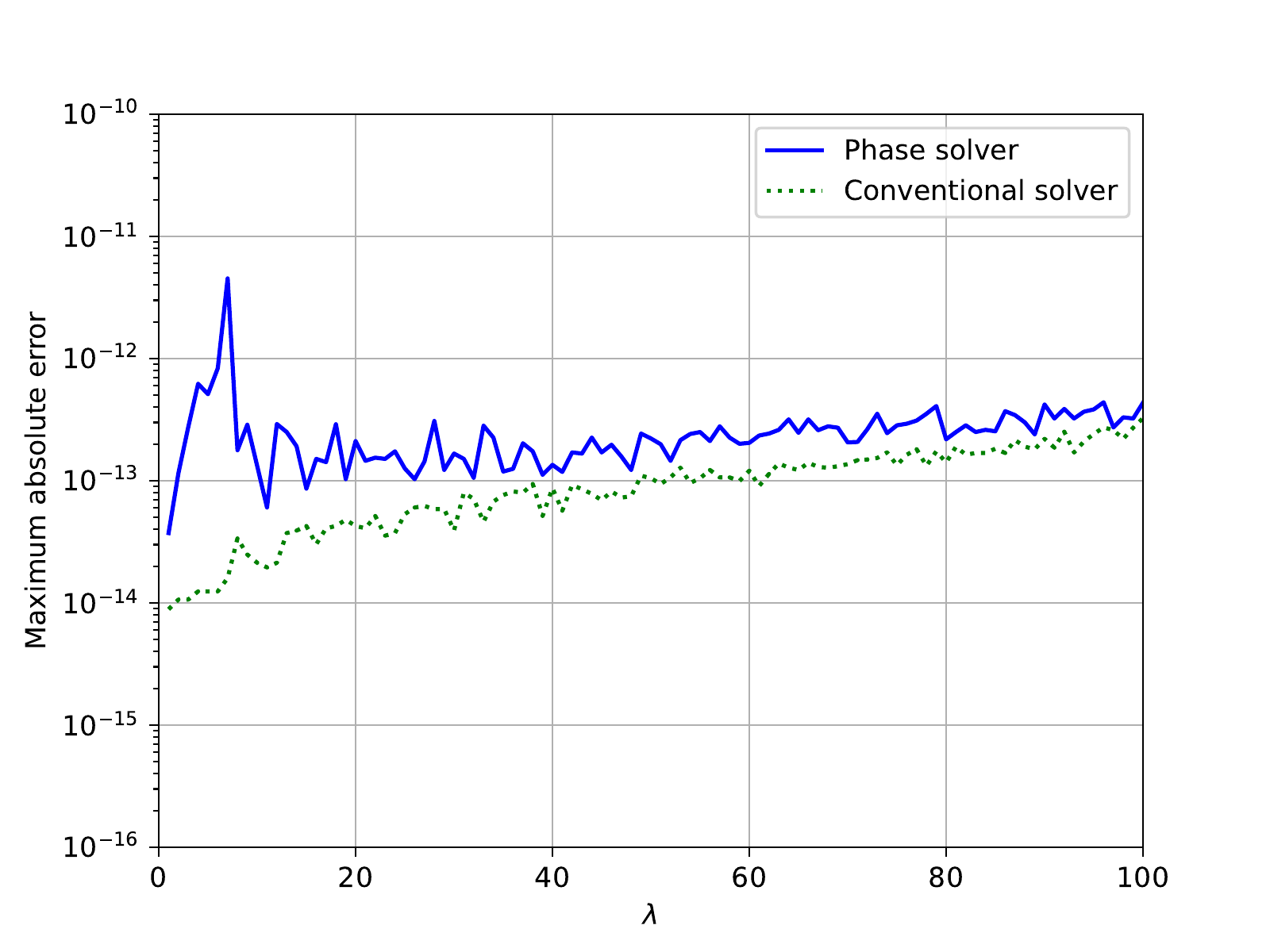}
\hfil
\includegraphics[width=.40\textwidth]{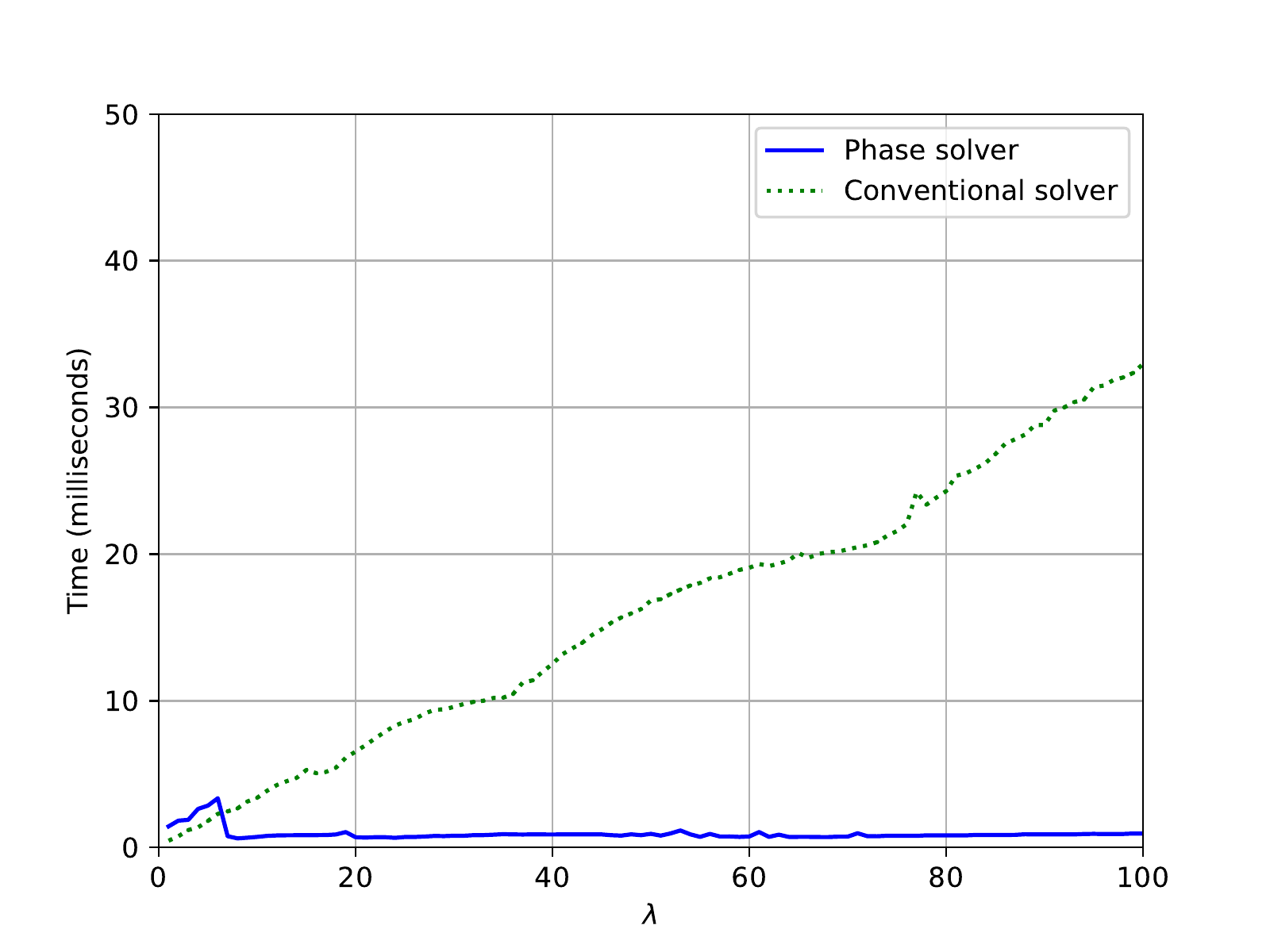}
\hfil

\hfil
\includegraphics[width=.40\textwidth]{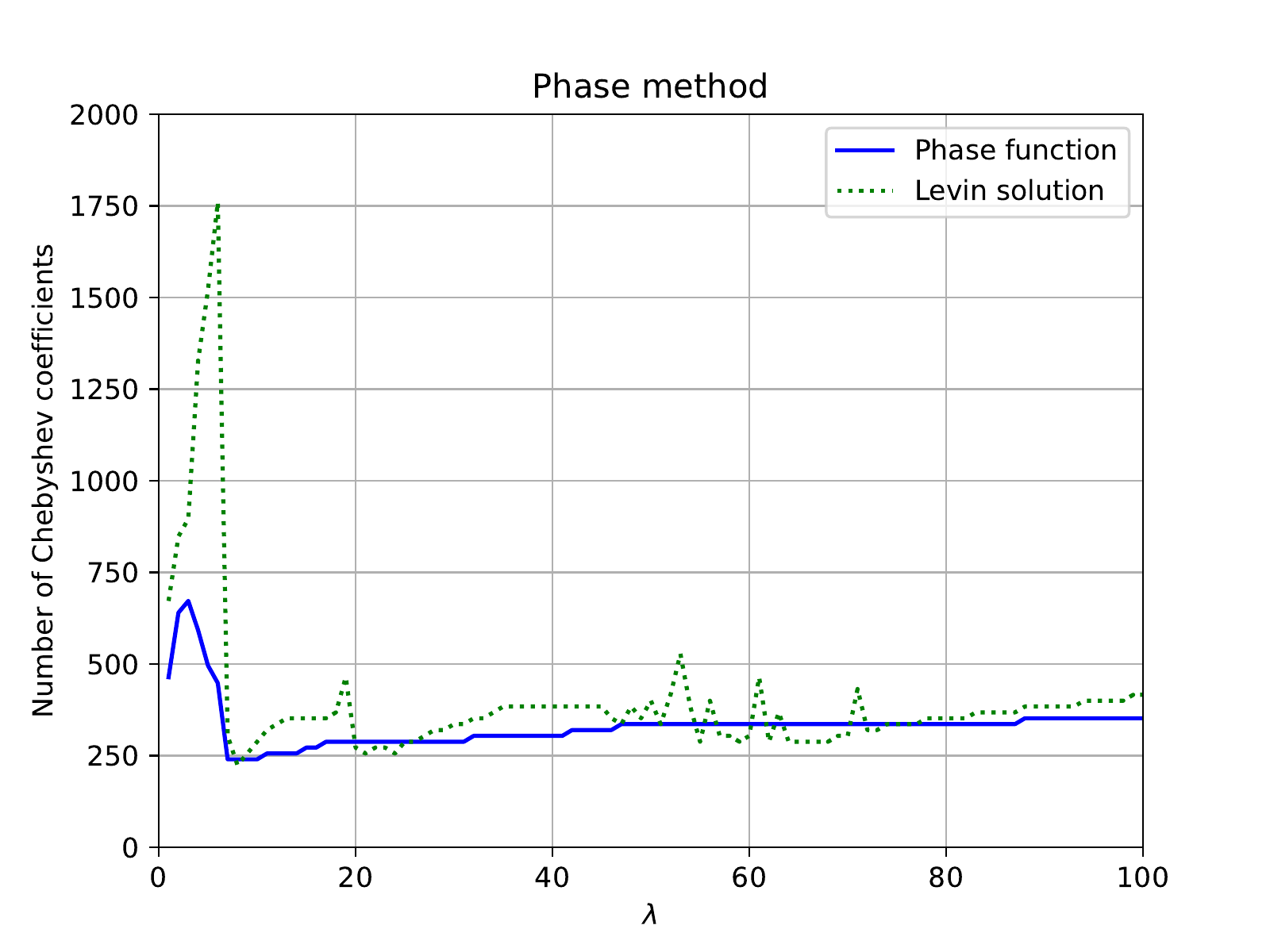}
\hfil
\includegraphics[width=.40\textwidth]{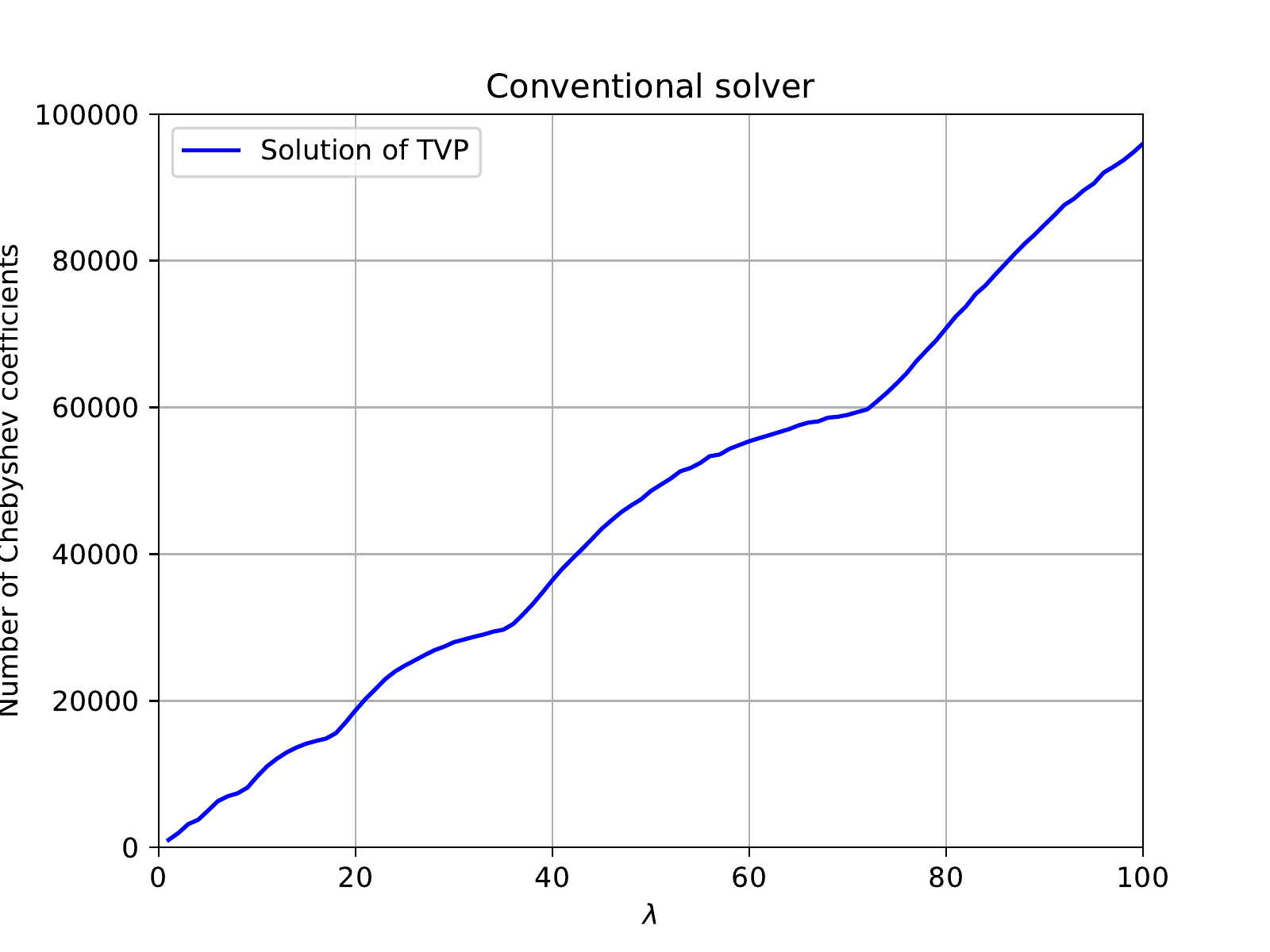}
\hfil

\caption{
The results of the second set of experiments of Subsection~\ref{section:experiments:airy},
which demonstrate the behavior of our algorithm for values of the parameter $\lambda$
between $1$ and $100$.
The plots in the first row compare the running time and accuracy
of the phase method with a conventional solver.  The plot in the lower left
give the number of Chebyshev coefficients used to represent the phase function
and the Levin solution, while the plot in the lower right gives the number of 
Chebyshev coefficients used to represent the solution produced by the conventional
solver.  
}
\label{figure:airyplots2}
\end{figure}

In our first set of experiments, we used both the phase function method and the conventional
approach described in \ref{section:algorithm:odesolver} to solve the terminal value problem
\begin{equation}
\left\{
\begin{aligned}
&y''(t) - \lambda^2 t  y(t) = \lambda^2 t^2, \ \ \ -10 < t< 0, \\
&y(0) = \frac{1}{3^{\frac{2}{3}} \Gamma\left(\frac{2}{3}\right)} \\
&y'(0) = -1 -\frac{\lambda^{\frac{2}{3}}\Gamma\left(-\frac{1}{3}\right)}{2\pi\ 3^{\frac{5}{6}}}.
\end{aligned}
\right.
\label{experiments:airy:ivp}
\end{equation}
It can be easily verified that the solution of (\ref{experiments:airy:ivp}) is
\begin{equation}
y(t) = -t + \mbox{Ai}\left(\lambda^{\frac{2}{3}} t\right),
\label{experiments:airy:sol}
\end{equation}
where $\mbox{Ai}$ denotes the Airy function of the first kind (see, for instance, \cite{Olver}
for a discussion of the Airy functions and a definition of the function $\mbox{Ai}$).

Because the condition number of the problem (\ref{experiments:airy:ivp})
and the condition number of evaluation of the function (\ref{experiments:airy:sol})
increase with $\lambda$, we took the precision parameter for the conventional
solver to be 
\begin{equation}
\epsilon_{\mbox{conventional}} = \max\left\{10^{-13}, \epsilon_0 \times \lambda\right\},
\end{equation}
where $\epsilon_0 \approx 2.220446049250313\times^{-16}$ is machine zero for IEEE double precision arithmetic.
As in all of the experiments discussed in this paper, we set the parameter specifying
the desired precision for the phase function and for the algorithm of Subsection~\ref{section:levin:algorithm}  to be
$\epsilon = 10^{-13}$.

\begin{figure}[t!!!!!!!!!!]
\centering
\hfil
\includegraphics[width=.40\textwidth]{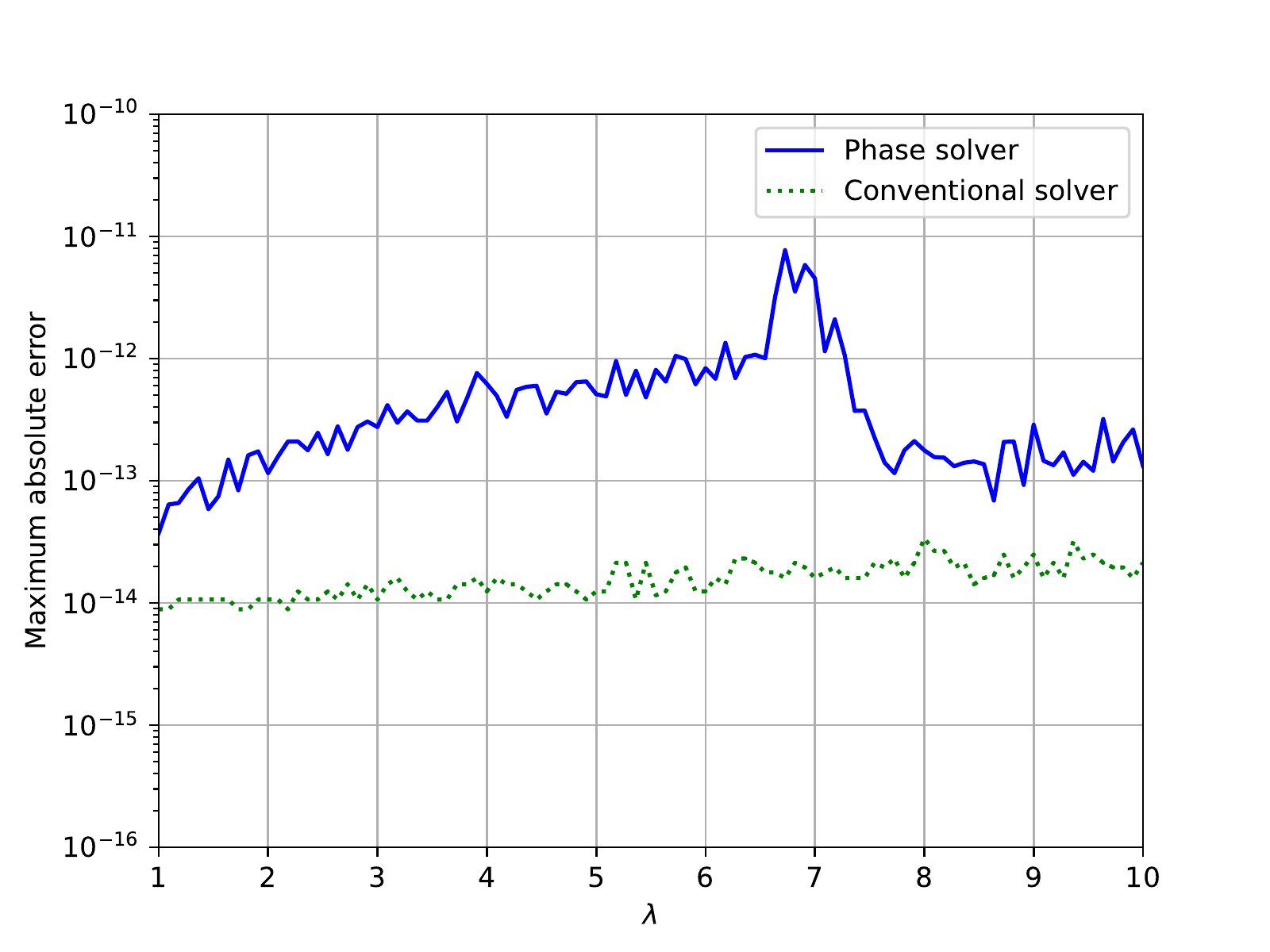}
\hfil
\includegraphics[width=.40\textwidth]{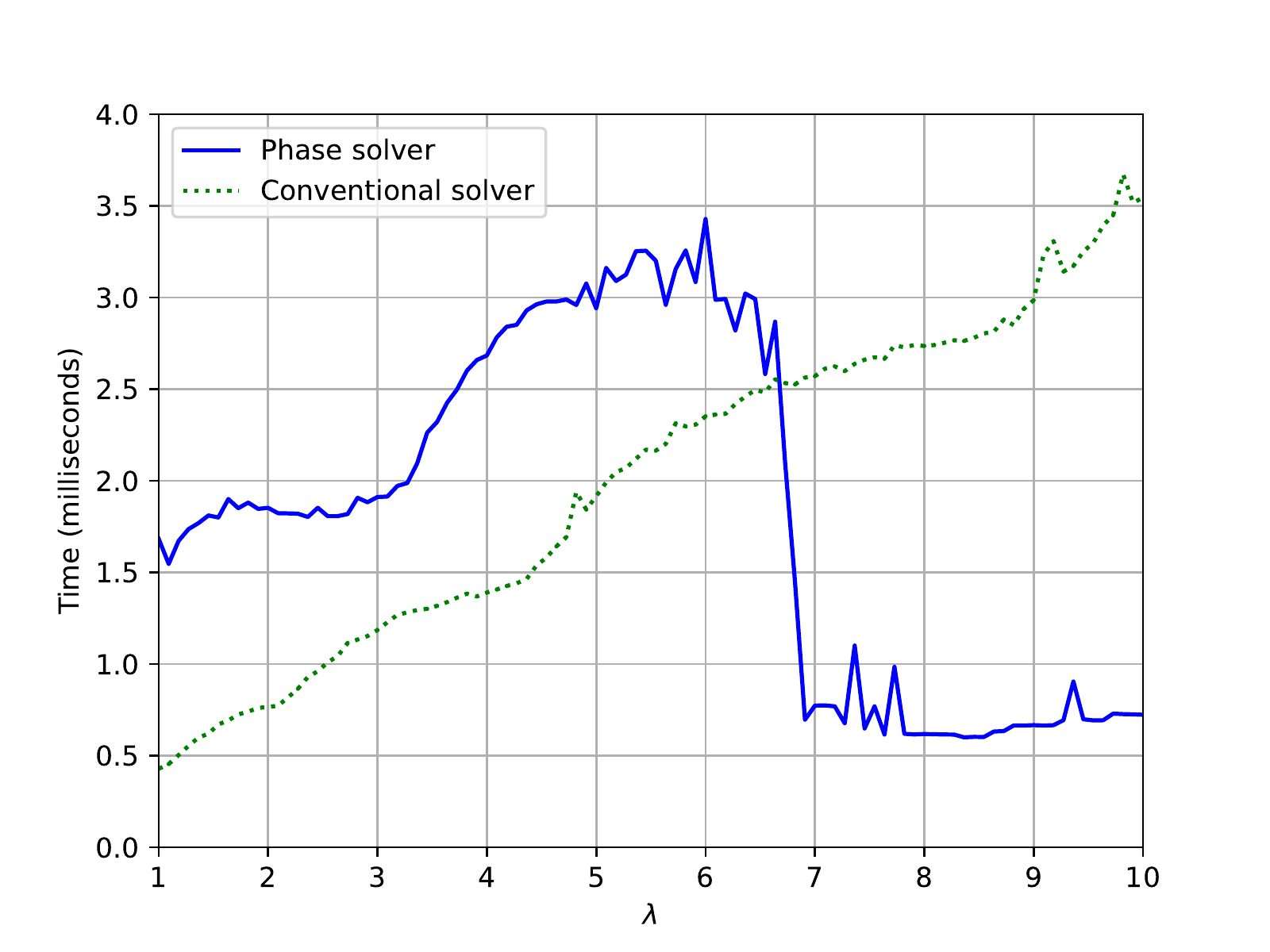}
\hfil

\hfil
\includegraphics[width=.40\textwidth]{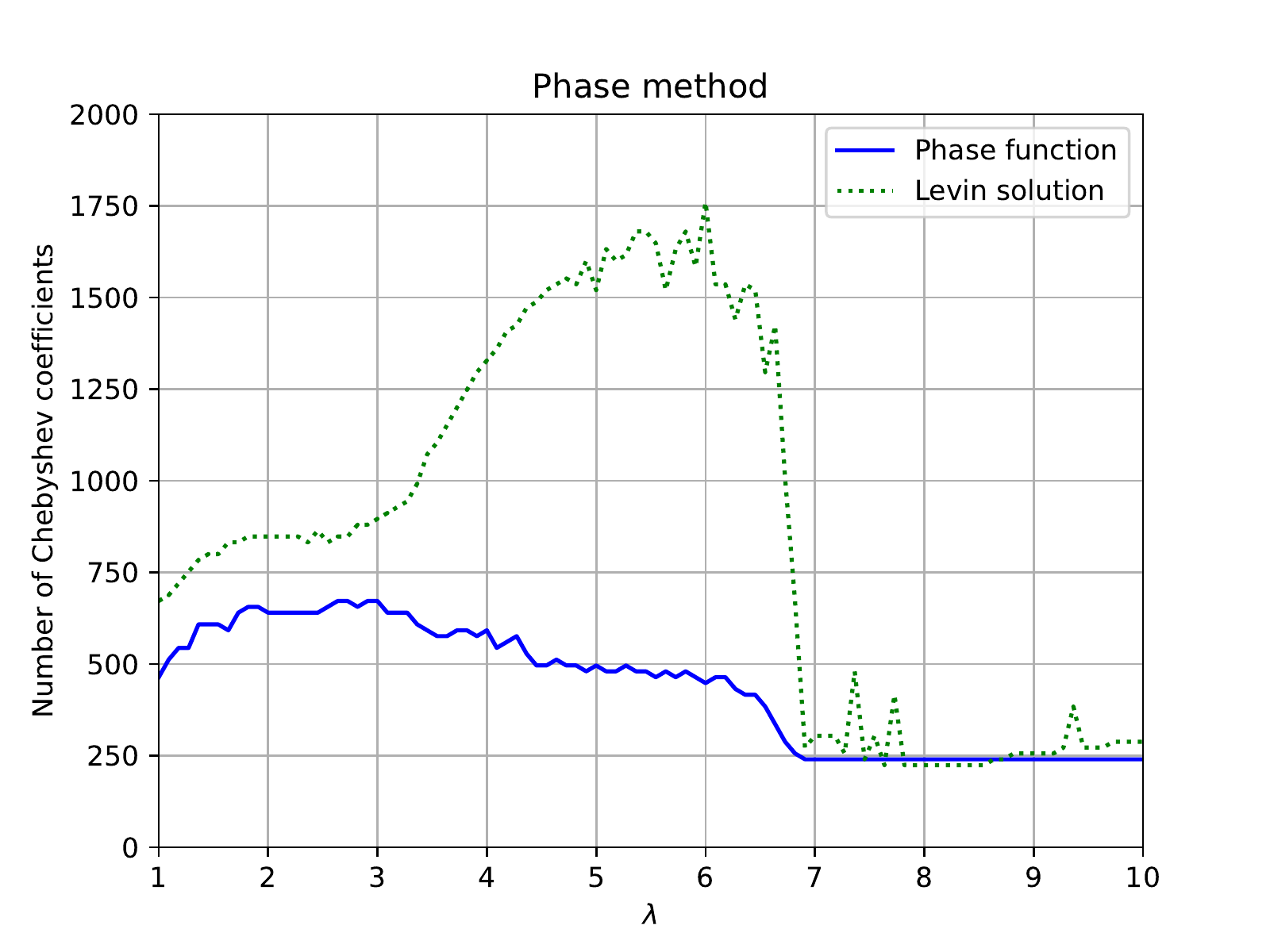}
\hfil
\includegraphics[width=.40\textwidth]{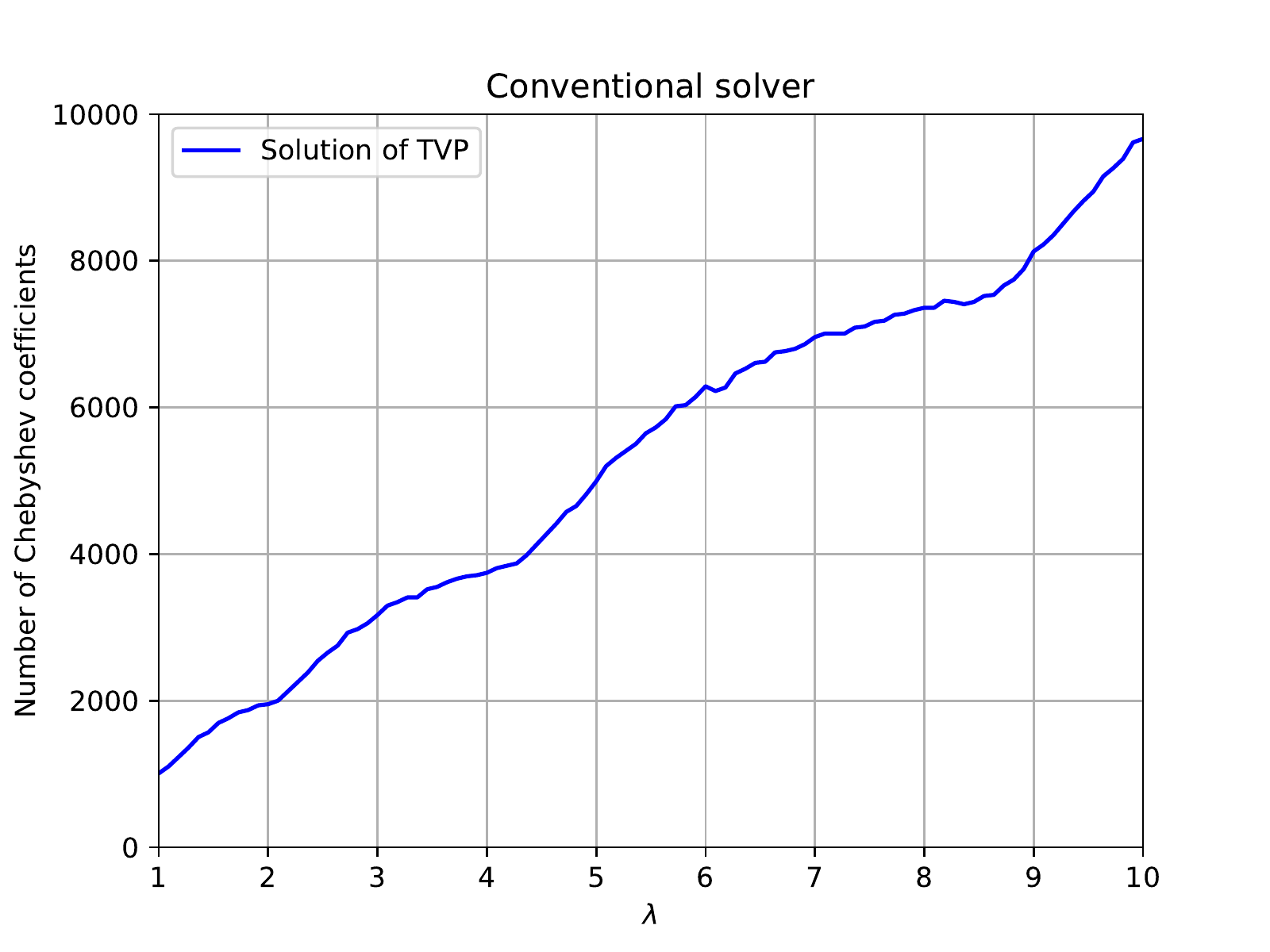}
\hfil

\caption{
The results of the third set of experiments of Subsection~\ref{section:experiments:airy},
which demonstrate the behavior of our algorithm for values of the parameter $\lambda$
between $1$ and $10$.
The plots in the first row compare the accuracy and running time 
of the phase method with a conventional solver.  The plot in the lower left
give the number of Chebyshev coefficients used to represent the phase function
and the Levin solution, while the plot in the lower right gives the number of 
Chebyshev coefficients used to represent the solution produced by the conventional
solver.  
}
\label{figure:airyplots3}
\end{figure}

We began our first experiment by sampling $m=100$ equispaced points $x_1,\ldots,x_m$ in the interval $[1,6]$.  Then, for each
$\lambda = 10^{x_1}, 10^{x_2}, \ldots, 10^{x_m}$, we solved (\ref{experiments:airy:ivp})
using both the phase function method and the conventional solver.
The wall clock time required by each method was measured, as was the number of coefficients in the piecewise
Chebyshev expansion used to represent the phase function and the total number of Chebyshev coefficients
in the expansions of the functions $p_1,\ldots,p_m$ produced by the algorithm of Subsection~\ref{section:levin:algorithm}.
In a mild abuse of terminology, we refer to this latter quantity in the rest of this section and in our figures
as the number of piecewise Chebyshev coefficients needed to represent the Levin solution.  Each of the two obtained solutions was then evaluated  
at $10,000$ equispaced points in the interval $(-10,0)$
and the largest observed absolute errors recorded.  The errors were, of course,
measured by comparison with the known solution (\ref{experiments:airy:sol}).  We measured absolute errors
rather than relative errors because (\ref{experiments:airy:sol}) is an oscillatory function with many zeros
on the interval $(-10,0)$.    The results are presented in Figure~\ref{figure:airyplots1}.  
The plot at the top of that figure gives the maximum observed absolute errors in the solutions
obtained by the phase function method and the conventional solver as functions of $\lambda$.
The plot on the left-hand side of the second row gives the time required by the phase method as
a function of $\lambda$, while the plot on the right-hand side of the second row gives the time
required by the conventional solver as a function of $\lambda$.
The plot on the left-hand side of the third row gives the number of coefficients
in the piecewise Chebyshev expansions of the phase function and the Levin solution.
The plot on the right-hand side of the third row
gives the number of coefficients in the piecewise Chebyshev expansion
of the solution of (\ref{experiments:airy:ivp}) produced by the conventional solver.

Two more experiments concerning the problem (\ref{experiments:airy:ivp}) were performed to test 
the behavior of the algorithm of this paper in the case of small values of $\lambda$.  
In the first of these experiments, we sampled
$m=100$ equispaced points $\lambda_1,\lambda_2, \ldots,\lambda_m$ in the interval $[1,100]$.
For each obtained value of $\lambda$, we repeated the procedure described above.
The results are given in Figure~\ref{figure:airyplots2}.
In the second of these experiments, we sampled
$m=100$ equispaced points $\lambda_1,\lambda_2, \ldots,\lambda_m$ in the interval $[1,10]$
and repeated the experiments described above for each obtained value of $\lambda$.
The results are shown in Figure~\ref{figure:airyplots3}.

\end{subsection}

\begin{subsection}{An initial value problem}
\label{section:experiments:two}

In our next experiment, we solved the initial value problem
\begin{equation}
\left\{
\begin{aligned}
&y''(t) + \left(\frac{\lambda^2}{0.01+t^2}\right) y(t) = \lambda^2 (1+t) \cos(13t^2), \ \ \ 0 <t < 1, \\
&y(0) = y'(0) = 1.
\end{aligned}
\right.
\label{experiments:two:problem}
\end{equation}
for various values of $\lambda$ using the phase function method.  More explicitly, we first
sampled $m=100$ points $x_1,x_2,\ldots,x_m$ in the interval $[1,6]$.  Then,
for each $\lambda=10^{x_1}, 10^{x_2}, \ldots, 10^{x_m}$, we used the algorithm
of this paper to solve (\ref{experiments:two:problem}).  For each $\lambda$, we recorded 
the wall clock time spent solving the problem and measured the error in the 
obtained  solution at $10,000$ equispaced points in the interval $[0,1]$.  The reference
solution was obtained by solving (\ref{experiments:two:problem}) using the conventional
solver running in quadruple precision (REAL*16) arithmetic.  
The results are presented in Figure~\ref{figure:exp2plots}.  

\begin{figure}[h!!!!!!!!!!!!!!!!!!!!!!]
\centering

\hfil
\includegraphics[width=.32\textwidth]{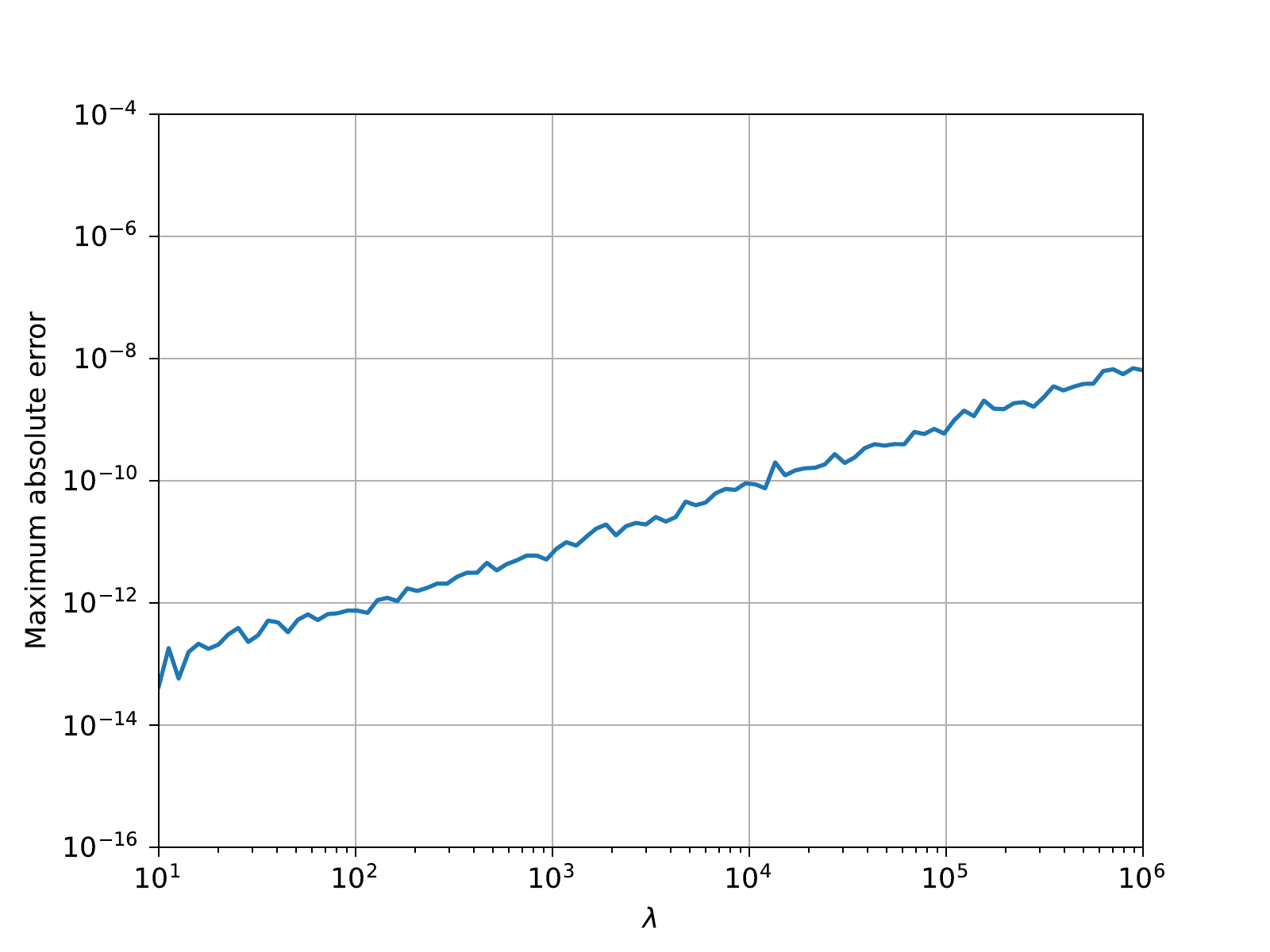}
\hfil
\includegraphics[width=.32\textwidth]{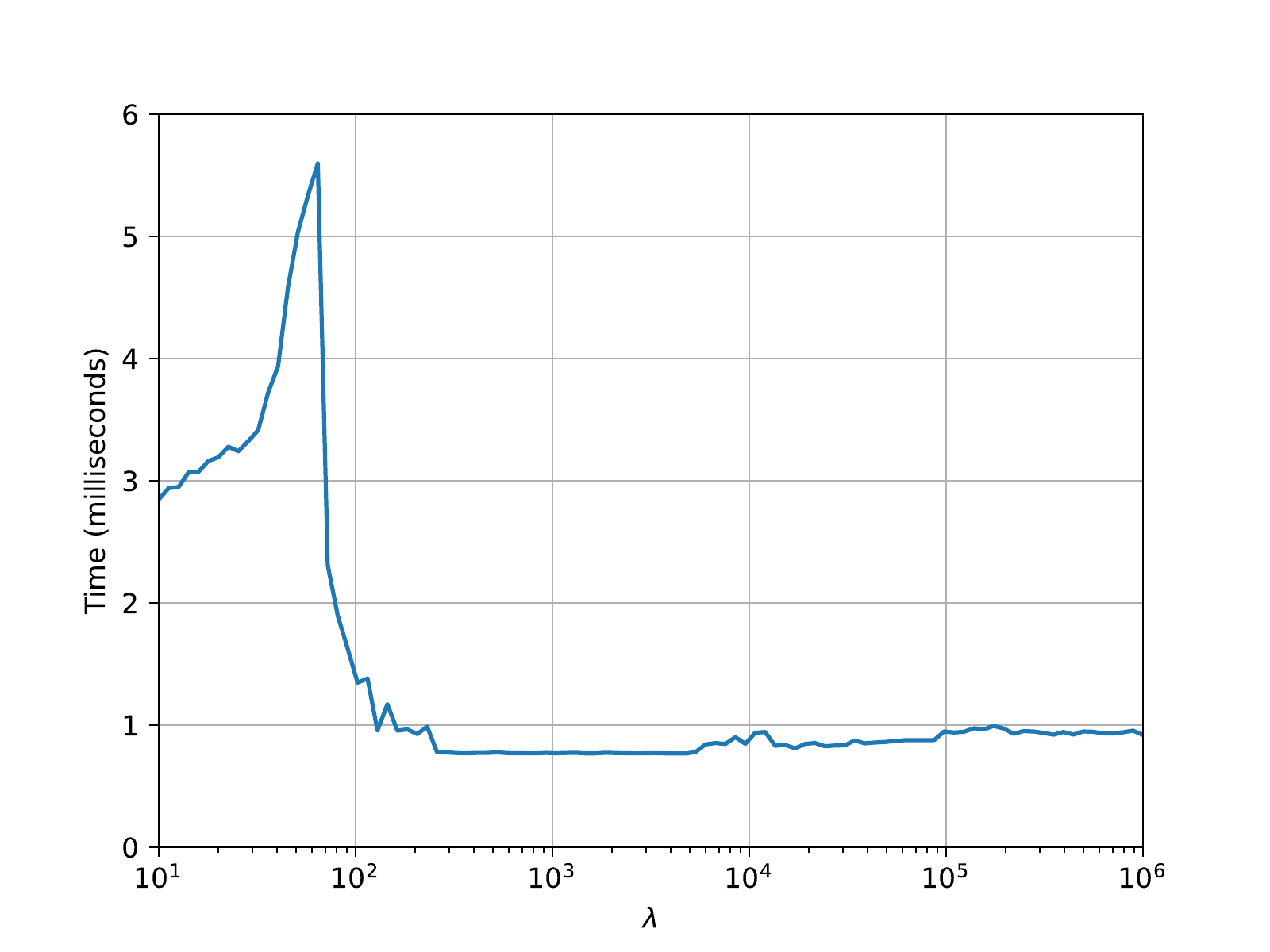}
\hfil
\includegraphics[width=.32\textwidth]{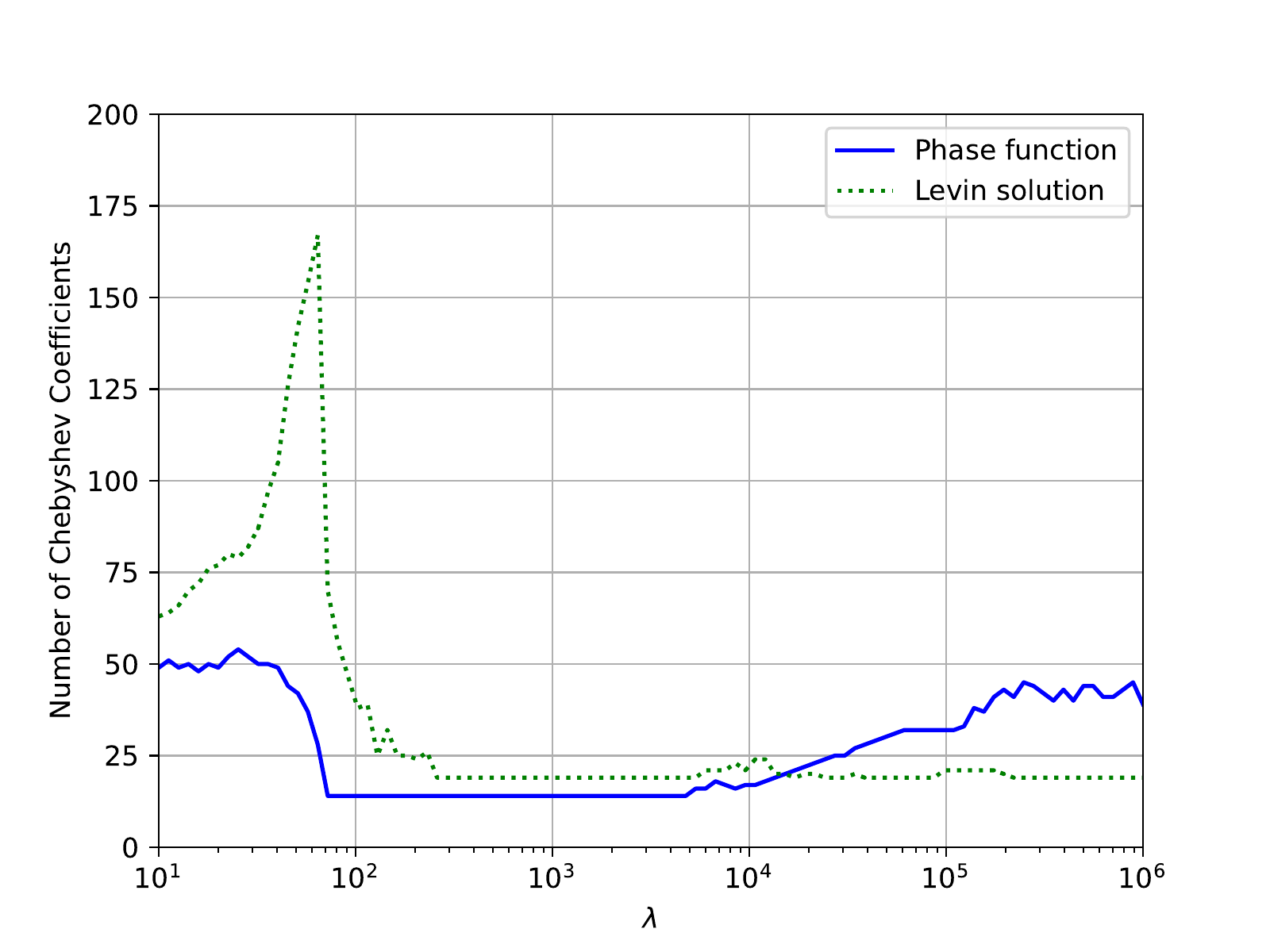}
\hfil

\caption{The results of the experiments of Subsection~\ref{section:experiments:two}.
The plot on the left
gives the maximum observed absolute error as a function of $\lambda$; the middle plot
gives the wall clock time required to solve (\ref{experiments:two:problem}) as a function
of $\lambda$; and the plot on the right-hand side gives the number of coefficients
in the piecewise Chebyshev expansions of the phase function and the Levin solution
as functions of the parameter $\lambda$.
}
\label{figure:exp2plots}
\end{figure}

\end{subsection}


\begin{subsection}{A boundary value problem}
\label{section:experiments:three}

In the experiment described in this subsection, we solved the boundary value problem
\begin{equation}
\left\{
\begin{aligned}
&y''(t) + \left(\frac{\lambda^3 \left(\frac{3}{2}+\cos\left(\log(\lambda)t\right)\right)}{1+\lambda\exp(t)} \right) y(t) = \frac{\lambda^2}{\sqrt{2+t}}, 
\ \ \ -1 <t < 1,\\
&y(-1) = y(1) = 0
\end{aligned}
\right.
\label{experiments:three:problem}
\end{equation}
for various values of $\lambda$.
More explicitly, we first
sampled $m=100$ points $x_1,x_2,\ldots,x_m$ in the interval $[1,6]$.  Then,
for each $\lambda=10^{x_1}, 10^{x_2}, \ldots, 10^{x_m}$, we used the algorithm
of this paper to solve (\ref{experiments:three:problem}).  For each $\lambda$, we recorded 
the wall clock time spent solving the problem and measured the error in the 
obtained  solution at $10,000$ equispaced points in the interval $[-1,1]$.  The reference
solution was obtained by solving (\ref{experiments:three:problem}) using the conventional
solver running in quadruple precision (REAL*16) arithmetic.  
The results are presented in Figure~\ref{figure:exp3plots}.  
We note that the coefficient in the differential equation appearing in
(\ref{experiments:three:problem}) becomes more oscillatory with $\lambda$, so we expect
the running time of our algorithm to grow with $\lambda$.  And while this is, in fact, the case,
the effect is remarkably mild.

\begin{figure}[h!]
\centering

\hfil
\includegraphics[width=.32\textwidth]{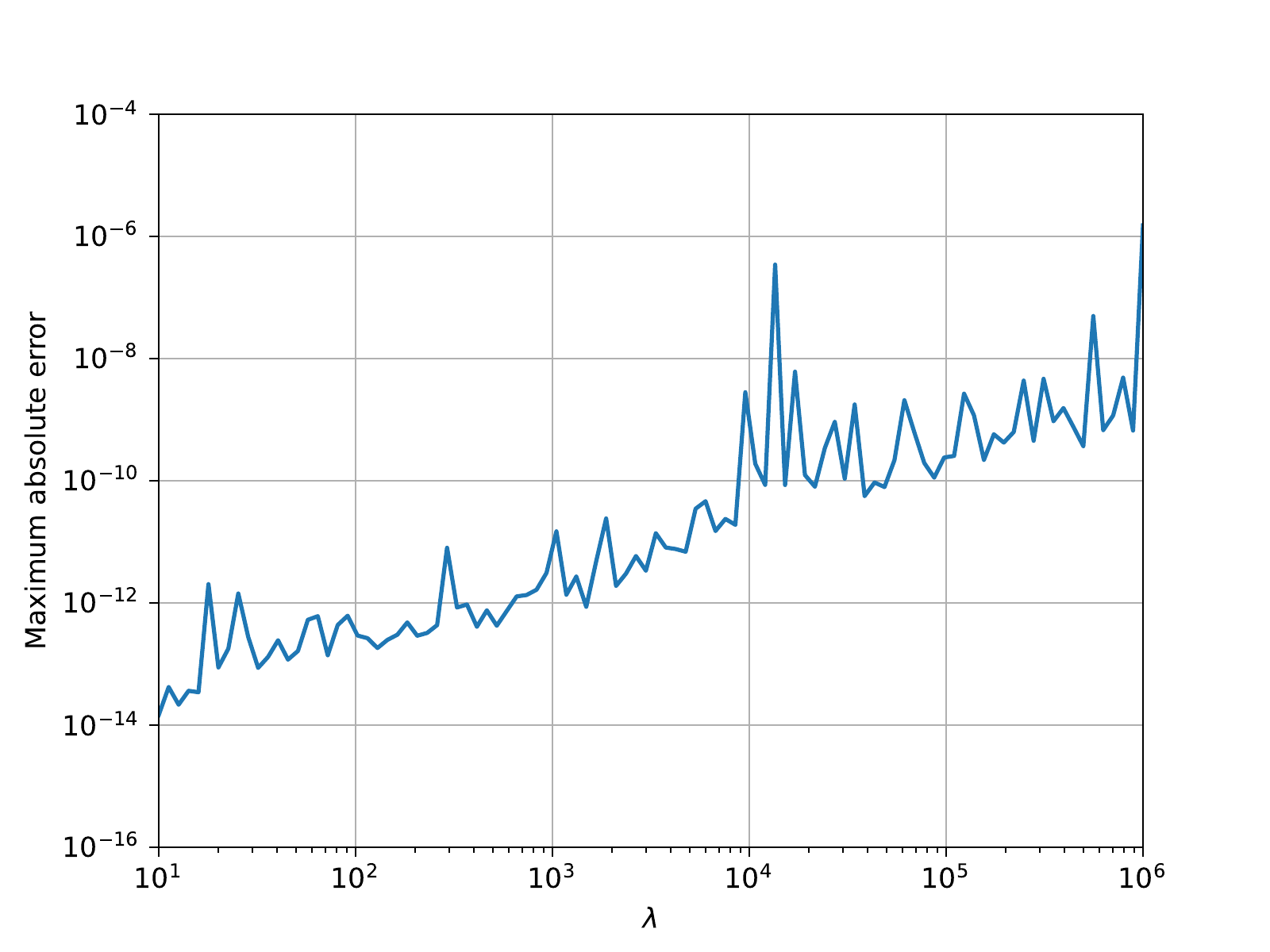}
\hfil
\includegraphics[width=.32\textwidth]{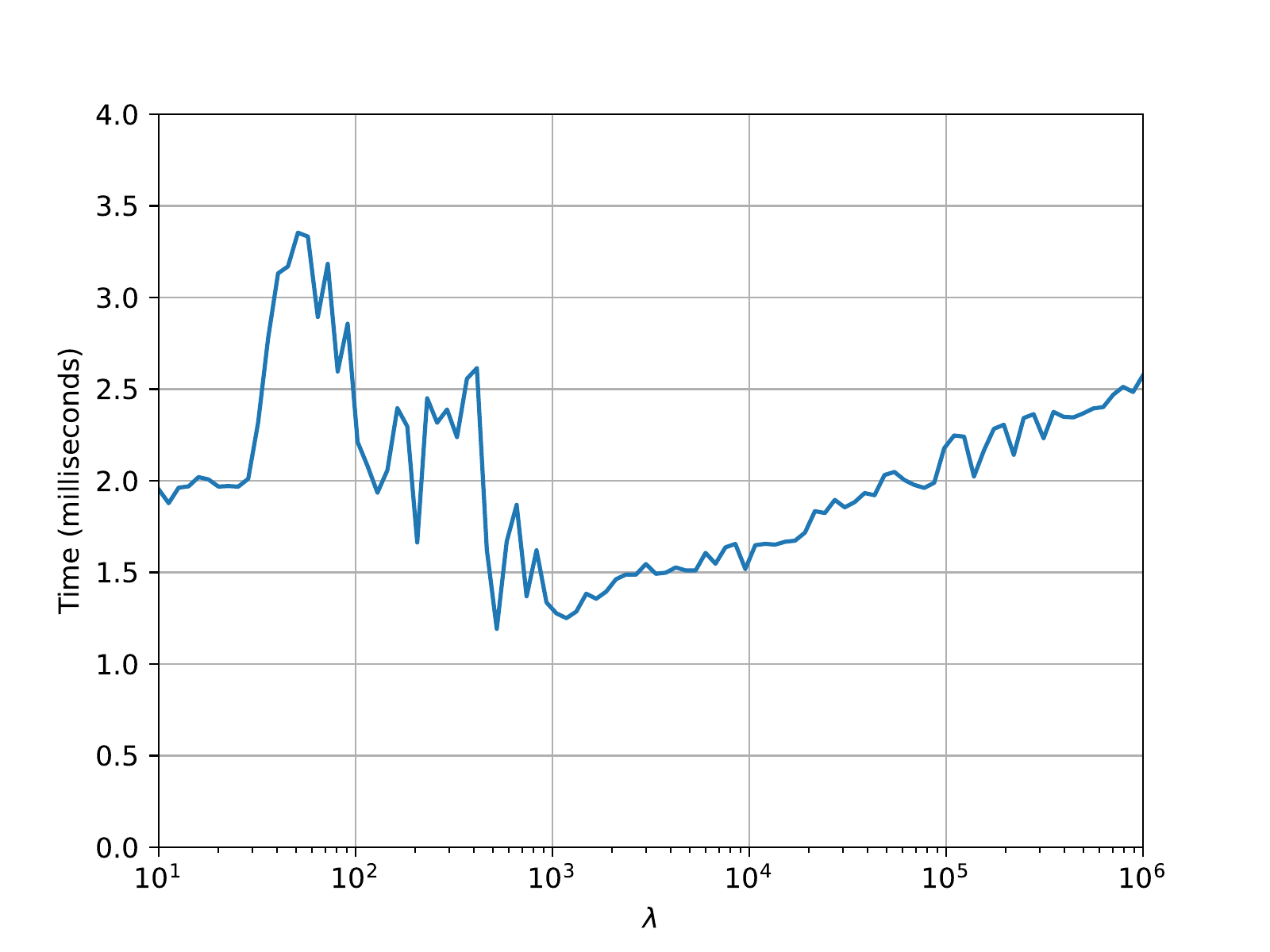}
\hfil
\includegraphics[width=.32\textwidth]{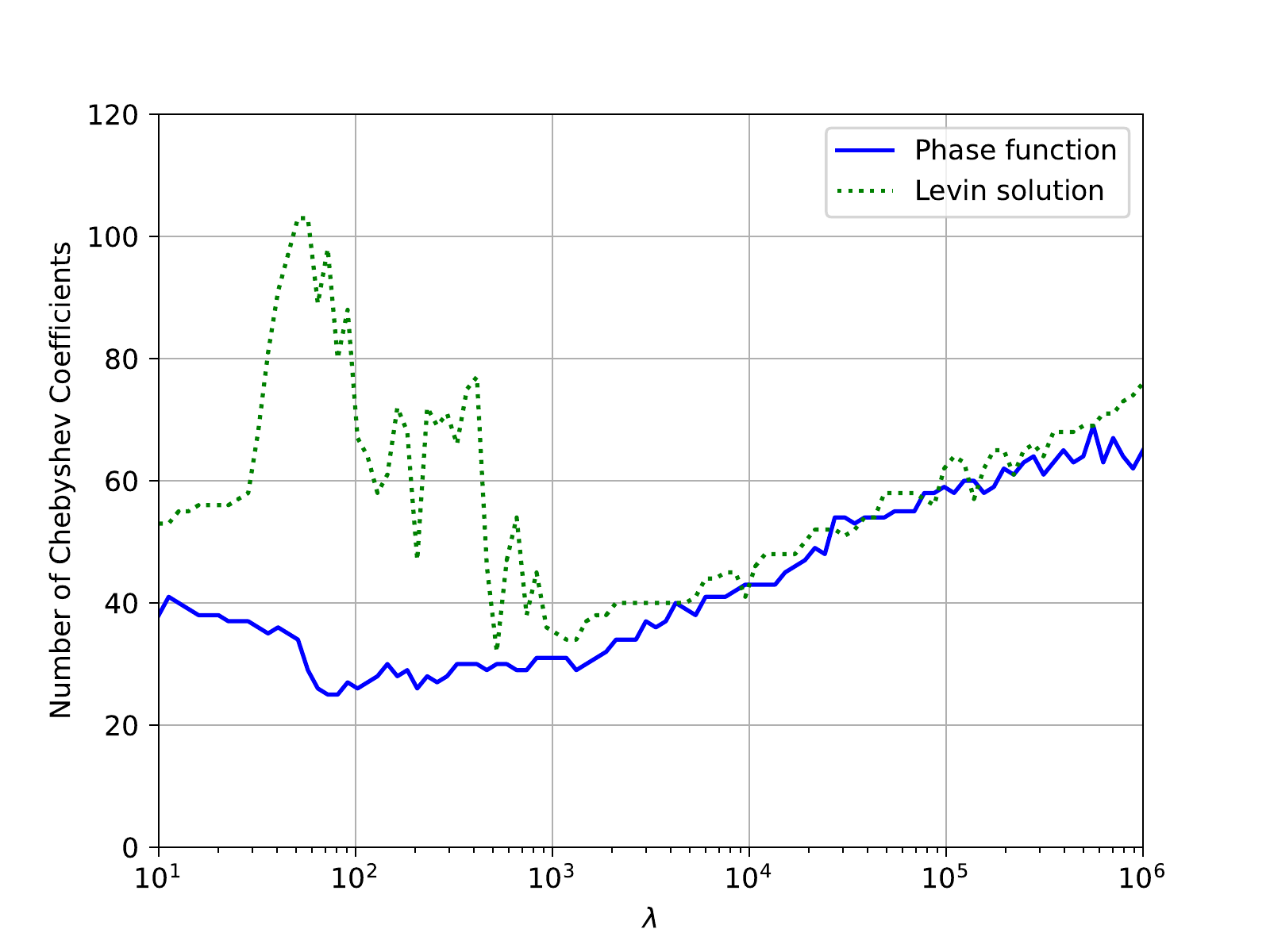}
\hfil

\caption{The results of the experiments of Subsection~\ref{section:experiments:three}.
The plot on the left
gives the maximum observed absolute error as a function of $\lambda$; the middle plot
gives the wall clock time required to solve (\ref{experiments:three:problem}) as a function
of $\lambda$; and the plot on the right-hand side gives the number of coefficients
in the piecewise Chebyshev expansions of the phase function and the Levin solution
as functions of the parameter $\lambda$.
}
\label{figure:exp3plots}
\end{figure}
\end{subsection}


\begin{subsection}{A boundary value problem with more complicated boundary conditions}
\label{section:experiments:four}

In a final experiment, we solved the boundary value problem
\begin{equation}
\left\{
\begin{aligned}
&y''(t) + \left( \frac{\lambda^2\left(2+t^2\cos(\lambda)\right)}{1+t^2}\right) y(t) = 
\lambda^2\cos(3t^2), \ \ \ -1 <t < 1,\\
&y(-1) = y(1) \\
&y'(-1) = y'(1)
\end{aligned}
\right.
\label{experiments:four:problem}
\end{equation}
for various values of $\lambda$.
We first sampled $m=100$ points $x_1,x_2,\ldots,x_m$ in the interval $[1,6]$.  Then,
for each $\lambda=10^{x_1}, 10^{x_2}, \ldots, 10^{x_m}$, we used the algorithm
of this paper to solve (\ref{experiments:four:problem}).  For each $\lambda$, we recorded 
the wall clock time spent solving the problem and measured the error in the 
obtained  solution at $10,000$ equispaced points in the interval $[-1,1]$.  The reference
solution was obtained by solving (\ref{experiments:four:problem}) using the conventional
solver running in quadruple precision (REAL*16) arithmetic.  
The results are presented in Figure~\ref{figure:exp4plots}.

\begin{figure}[h!]
\centering

\hfil
\includegraphics[width=.32\textwidth]{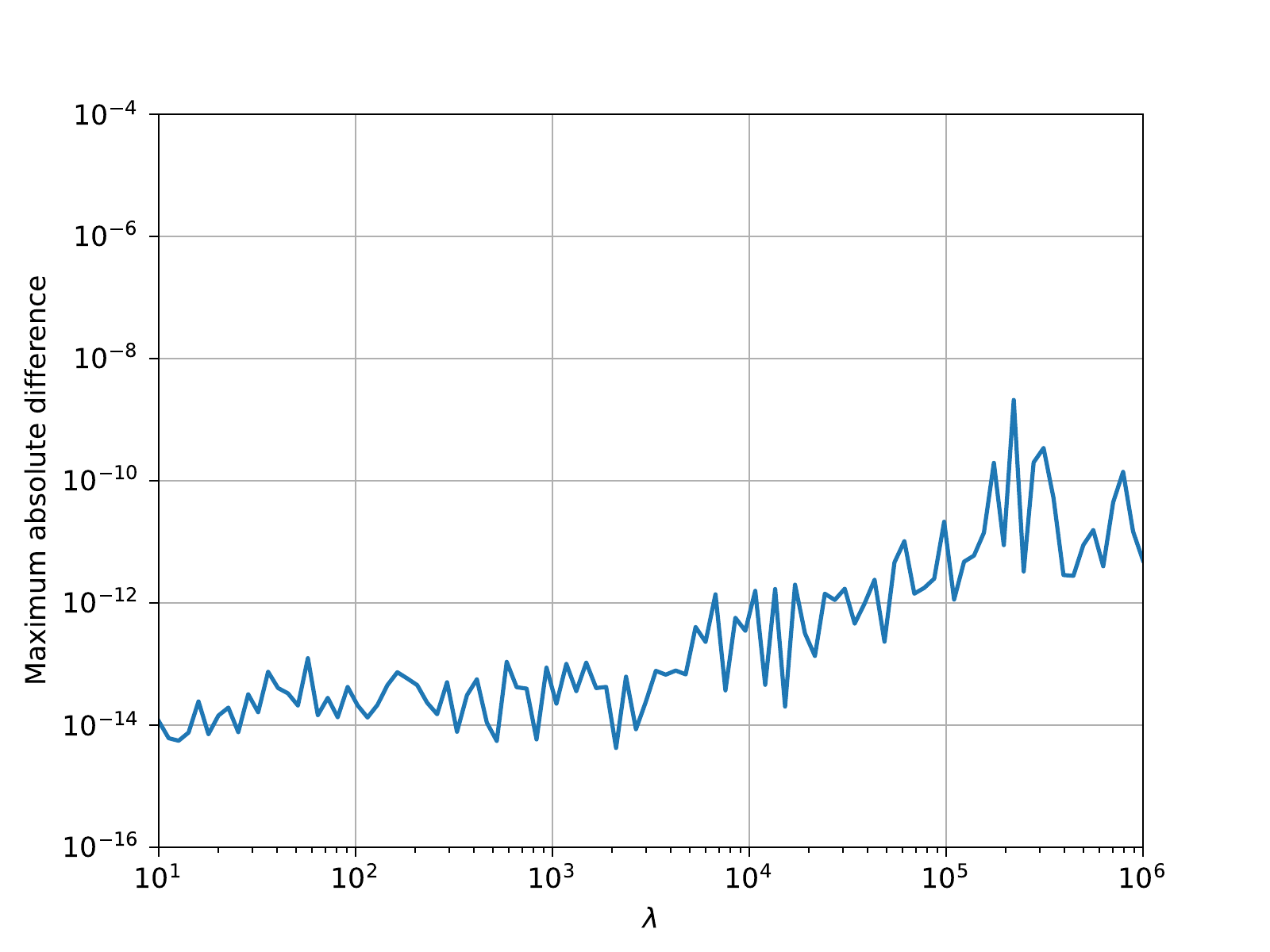}
\hfil
\includegraphics[width=.32\textwidth]{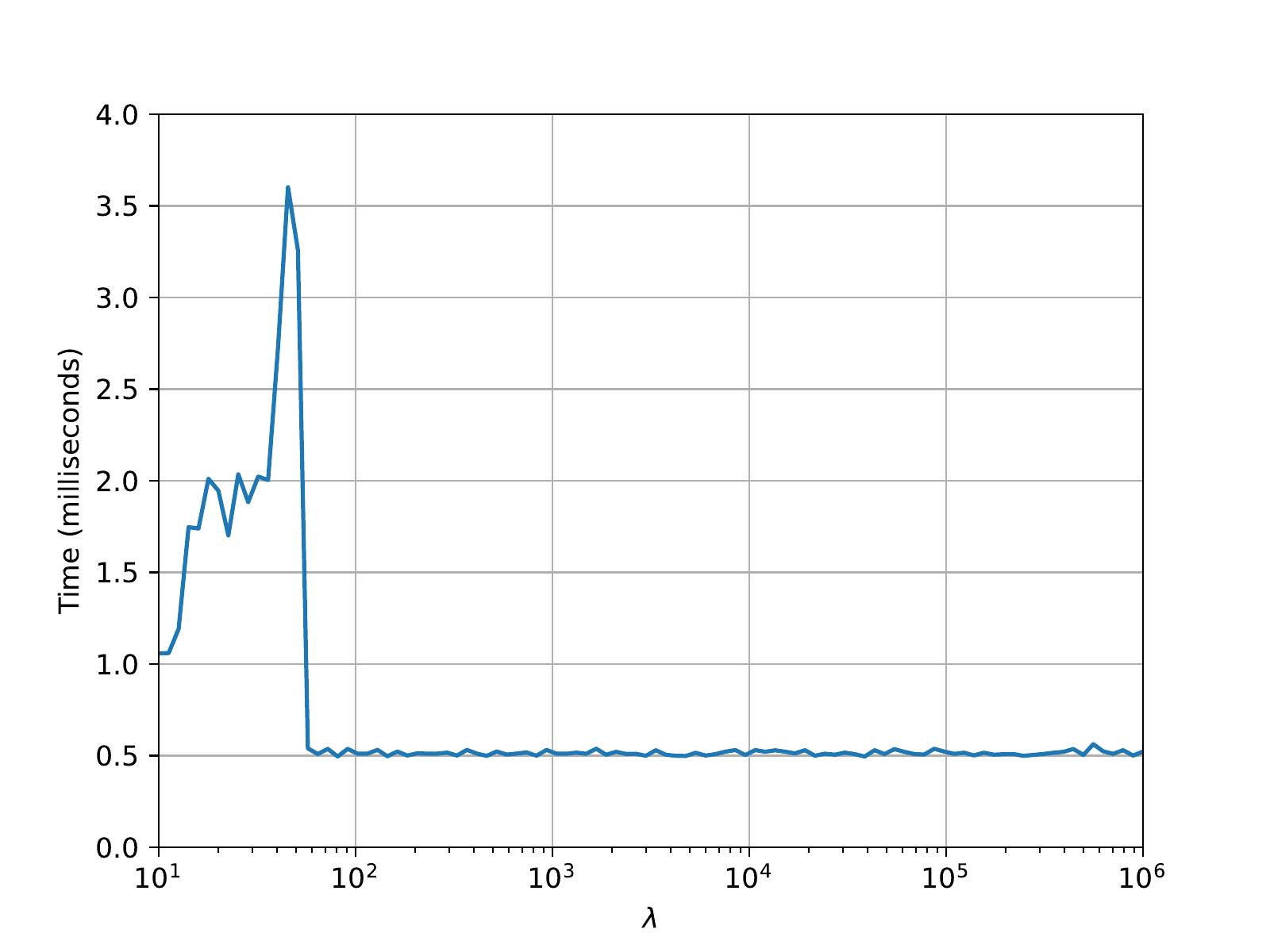}
\hfil
\includegraphics[width=.32\textwidth]{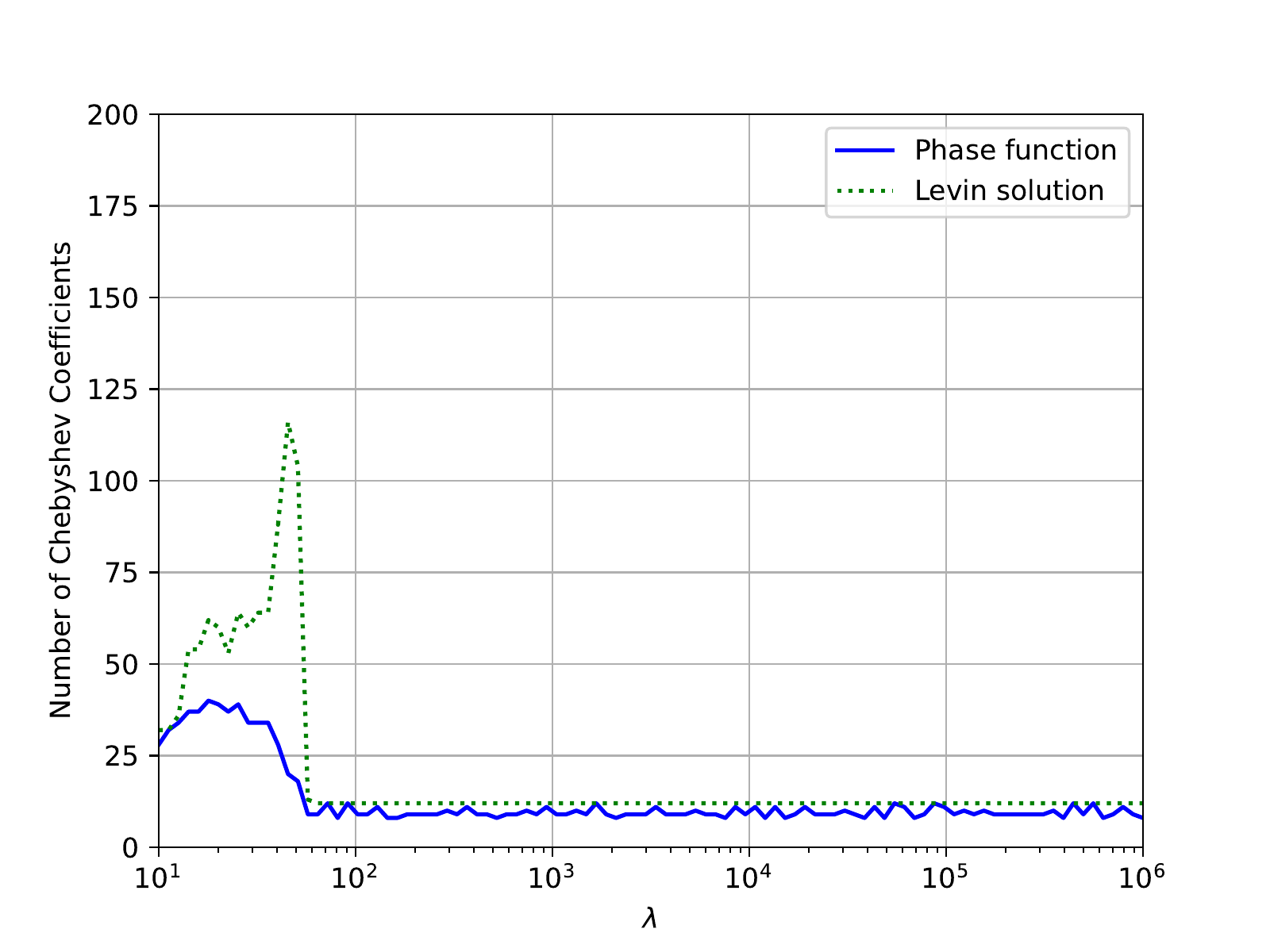}
\hfil

\caption{The results of the experiments of Subsection~\ref{section:experiments:four}.
The plot on the left
gives the maximum observed absolute error as a function of $\lambda$; the middle plot
gives the wall clock time required to solve (\ref{experiments:four:problem}) as a function
of $\lambda$; and the plot on the right-hand side gives the number of coefficients
in the piecewise Chebyshev expansions of the phase function and the Levin solution
as functions of the parameter $\lambda$.
}
\label{figure:exp4plots}
\end{figure}

\end{subsection}

\end{section}

\begin{section}{Conclusions}
\label{section:conclusions}

We have introduced a scheme for solving second order inhomogeneous linear ordinary differential equations of the form
\begin{equation}
y''(t) + q(t) y(t) = f(t),\ \ \ \ \ a < t < b,
\label{conclusion:helmholtz}
\end{equation}
in the case in which $f$ is real-valued and slowly varying and $q$ is positive and  slowly varying.
Unlike standard solvers, its running time is 
independent of the magnitude of $q$ and, unlike asymptotic methods, the accuracy it obtains is consistent
with the condition number of the problem even when  $q$ is of small magnitude.

We chose the approach of this paper over another, more direct, Levin-type method
precisely because we wanted an algorithm which achieves high-accuracy regardless
of the magnitude of $q$. When $f$ and $q$ are slowly varying and $q$ is of sufficiently large magnitude, the 
equation (\ref{introduction:inhom})
admits a nonoscillatory solution and applying an algorithm analogous to that of 
Subsection~\ref{section:phase:algorithm} to it
will result in a piecewise Chebyshev expansion representing that solution.
However, such an approach will fail when $q$ is of small magnitude for
the reasons discussed at the end of Subsection~\ref{section:phase:algorithm}.
Nonetheless, it would be of some interest to determine when this more direct
approach is preferable to the algorithm of this paper.

As mentioned in the introduction, the representation (\ref{introduction:general})
of the general form of the solution of (\ref{introduction:inhom})
is numerically unstable when $q$ is negative and of large magnitude
on some or all of the interval $[a,b]$.  In the authors' view, it is unlikely
that an algorithm of the type described this paper can be successfully
applied to such problems in general.
 Nonetheless, it is most likely possible 
to further develop the method of this article so that it applies to important classes
of physically relevant problems for which the coefficient $q$ is negative on
some or all of the interval $[a,b]$ and this is the subject of ongoing work by the authors.

\end{section}

\begin{section}{Acknowledgements}
KS was supported in part by the NSERC Discovery 
Grants RGPIN-2020-06022 and DGECR-2020-00356.  
JB was supported in part by NSERC Discovery grant  RGPIN-2021-02613.
\end{section}

\bibliographystyle{acm}
\bibliography{odesolve.bib}

\appendix
\begin{section}{An adaptive Chebyshev spectral solver for ordinary differential equations}
\label{section:algorithm:odesolver}

The algorithm of this paper entails solving several ordinary differential equations.
We use a fairly standard adaptive Chebyshev spectral solver to do so.  We now briefly describe its operation
in the case of the initial value problem
\begin{equation}
\left\{
\begin{aligned}
\bm{y}'(t) &= F(t,\bm{y}(t)), \ \ \ a < t < b,\\
\bm{y}(a) &= \bm{v}
\end{aligned}
\right.
\label{algorithm:system}
\end{equation}
where $F:\mathbb{R}^{n+1} \to \mathbb{R}^n$ is smooth and $\bm{v} \in \mathbb{R}^n$.
The solver can be easily modified to apply to a terminal value problem.

The solver takes as input a positive integer $k$, a tolerance parameter $\epsilon$, an interval $(a,b)$, 
the vector $\bm{v}$ and a 
 subroutine for evaluating the function $F$.  It outputs $n$ piecewise $k^{th}$ order Chebyshev expansions,
one for each of the components $y_i(t)$ of the solution $\bm{y}$ of (\ref{algorithm:system}).

The solver maintains two lists of subintervals of $(a,b)$: one consisting of what we term ``accepted subintervals''
and the other of subintervals which have yet to processed.  A subinterval is accepted if the solution
is deemed to be adequately represented by a $k^{th}$ order Chebyhev expansion on that subinterval.
Initially, the list of accepted subintervals is empty and the list of 
subintervals to process contains the single interval $(a,b)$
It then operates as follows until the list of subintervals to process is empty:
\begin{enumerate}

\item
Find, in the list of subinterval to process, the interval $(c,d)$ such that
$c$ is as small as possible and remove this subinterval from the list.

\item
Solve the initial value problem
\begin{equation}
\left\{
\begin{aligned}
\bm{u}'(t) &= F(t,\bm{u}(t)), \ \ \ c< t < d,\\
\bm{u}(c) &= \bm{w}
\end{aligned}
\right.
\label{algorithm:ivp2}
\end{equation}
If $(c,d) = (a,b)$, then we take $\bm{w}=\bm{v}$.  Otherwise,
the value of the solution at the point $c$ has already been approximated, and we use that estimate
for $\bm{w}$ in (\ref{algorithm:ivp2}).

If the problem is linear, a straightforward integral equation method (see, for instance, \cite{GreengardSolver})
is used to solve (\ref{algorithm:ivp2}).  Otherwise, 
the trapezoidal method (see, for instance, \cite{Ascher}) is first used to produce an initial
approximation $\bm{y_0}$ of the solution and then Newton's method is applied to refine it.
The linearized problems are solved using an integral equation method.

In any event, the result is a set of $k^{th}$ order Chebyshev expansions 
\begin{equation}
u_i(t)  \approx \sum_{j=0}^k \lambda_{ij}\ T_j\left(\frac{2}{d-c} t + \frac{c+d}{c-d}\right),\ \ \ i=1,\ldots,n,
\label{algorithm:exps}
\end{equation}
approximating  the components $u_1,\ldots,u_n$ of the solution of (\ref{algorithm:ivp2}).

\item
Compute the quantities
\begin{equation}
\frac{\sqrt{\sum_{j=\lfloor k/2 \rfloor+1}^k \lambda_{ij}^2}}{\sqrt{\sum_{j=0}^k \lambda_{ij}^2}}, \ \ \ i=1,\ldots,n,
\end{equation}
where the $\lambda_{ij}$ are the coefficients in the expansions (\ref{algorithm:exps}).
If any of the resulting values is  larger than $\epsilon$,
then we split the subinterval into two halves $\left(c,\frac{c+d}{2}\right)$ and 
$\left(\frac{c+d}{2},d\right)$ and place them on the list of subintervals to process.  Otherwise, we place the subinterval
$(c,d)$ on the list of accepted subintervals.

\end{enumerate} 

At the conclusion of this procedure,  we have $k^{th}$ order Chebyshev expansions
for each component of the solution, with the list of accepted subintervals determining the
partition for each expansion.

\end{section}

\end{document}